%% file: CountWords-arxiv-2019-10-31.tex
\documentclass{article}
\usepackage{amsmath,amssymb,amsthm,amsfonts,mathrsfs}
\usepackage{geometry}
\usepackage{enumerate}
\usepackage{xcolor}
\usepackage{hyperref}
\hypersetup{
    colorlinks=true,
    linkcolor=blue,
    filecolor=blue,      
    urlcolor=blue,
    citecolor=blue!20!black,
}
 
\usepackage{tikz}
\usepackage{tikz-cd}
\usetikzlibrary{calc, intersections}
\usetikzlibrary{arrows}
\usetikzlibrary{matrix}
\usetikzlibrary{chains}
\makeindex
\usepackage{fancyhdr}

\DeclareFontEncoding{LS1}{}{}
\DeclareFontSubstitution{LS1}{stix}{m}{n}
\DeclareSymbolFont{symbols2}{LS1}{stixfrak} {m} {n}
\DeclareMathSymbol{\operp}{\mathbin}{symbols2}{"A8}

\newcommand{\iprod}[2]{[#1\, ,#2]}

\newcommand{\scalar}[2]{\langle #1\, ,#2\rangle}

\newcommand{\UGn}{U(\ma{G}_n)}
\DeclareMathOperator{\range}{Rg}

\newcommand{\define}[1]{\emph{#1}}

\definecolor{myred1}{rgb}{0.55,0.00,0.00}
\definecolor{myred2}{rgb}{1,0.06,0.06}
\definecolor{myblue}{rgb}{0.00,0.00,0.55}

\newcommand{\m}[1]{{\color{blue}\bfseries #1}}

\newtheorem{theorem}{Theorem}[section]
\newtheorem{prop}[theorem]{Proposition}
\newtheorem{lemma}[theorem]{Lemma}
\newtheorem{cor}[theorem]{Corollary}
\theoremstyle{definition}
\newtheorem{definition}[theorem]{Definition}
\newtheorem{example}[theorem]{Example}

\theoremstyle{remark}
\newtheorem{remark}[theorem]{Remark}
\numberwithin{equation}{section}

\newcommand{\op}{\operatorname}
\newcommand{\ma}{\mathcal}
\newcommand{\ben}{\begin{equation}}
\newcommand{\een}{\end{equation}}
\newcommand{\bena}{\begin{equation*}}
\newcommand{\eena}{\end{equation*}}
\newcommand{\To}{\longrightarrow}

\def\ZZ{\mathbb{Z}}
\def\NN{\mathbb{N}}

\def\FF{\mathbb{F}}
\newcommand{\mG}{\mathcal G}
\newcommand{\mL}{\mathcal L}
\newcommand{\mP}{\mathcal P}

\DeclareMathOperator{\linspan}{span}
\DeclareMathOperator{\supp}{supp}
\DeclareMathOperator{\wt}{wt}

\textwidth 16cm \textheight 21cm \fontsize{24}{30}

\makeindex

\begin{document}

\title{Calculating the dimension of the universal embedding of the symplectic dual polar space using languages}

\author{
C. Segovia\footnote{
Instituto de Matem\'aticas UNAM-Oaxaca,
csegovia@matem.unam.mx}\, and\,
M. Winklmeier\footnote{
Departamento de Matem\'aticas, Universidad de los Andes, Bogot\'a, Colombia,
mwinklme@uniandes.edu.co}
}

\maketitle

\begin{abstract}
   The main result of this paper is the construction of a bijection of the set of words in so-called standard order of length $n$ formed by four different letters and the set $\ma{N}^n$ of all subspaces of a fixed $n$-dimensional maximal isotropic subspace of the $2n$-dimensional symplectic space $V$ over $\FF_2$ which are not maximal in a certain sense.
   Since the number of different words in standard order is known, this gives an alternative proof for the formula of the dimension of the universal embedding of a symplectic dual polar space $\mathcal G_n$.
   Along the way, we give formulas for the number of all $n$- and $(n-1)$-dimensional totally isotropic subspaces of $V$.
\end{abstract}

\maketitle

\section*{Introduction}  
\label{intro}
Configurations of points and lines are of significant importance since they occur for instance
as designs in combinatorics, geometry and algebra.  
These structures have been extensively documented in \cite{Levi,Branko,PS}, and historically, projective geometry has provided important examples like the Fano plane \cite{Dem}. 
The configurations induced by a projective geometry are completely characterized by a set of axioms for its points and lines, and 
analogously we can find the configurations induced by the (dual) polar spaces. 
The axiomatic formulation of a polar space was given in \cite{francis}, 
while the axioms for a dual polar space were developed in \cite{cameron}.
An important example of a polar space is the set of all totally isotropic subspaces of a given symplectic space whereas the set of all maximal isotropic subspaces form a dual polar space. 

In this article we consider a symplectic space $V$ of dimension $2n$. 
We denote by $\ma{P}_n$ the set of all maximal totally isotropic subspaces of $V$ and by $\ma{L}_n$ the set of all totally isotropic subspaces of dimension $n-1$. 
They form a configuration of points and lines $\ma{G}_n=(\ma{P}_n,\ma{L}_n)$ called the \define{symplectic dual polar space}, where the incidence relation is given by inclusion of the subspaces.
In the case when $V$ is a $\FF_2$-vector space, this structure is completely understood and there is a vast literature on this matter \cite{BCN,bart,diageo}. 
The case $n=2$ is of great importance because it gives the self-dual configuration called Cremona-Richmond configuration \cite{crem,rich} whose exciting history can be found in \cite{Bak1,Bak2}. 
In Figure \ref{f1} we show the Cremona-Richmond configuration, which has fifteen points and fifteen lines and such that every point is contained in exactly three lines and every line contains exactly three different points.
\begin{figure}\label{f1}
   \begin{center}
      \input{CR.tex}
   \end{center}
   \caption{The Cremona-Richmond configuration.}
\end{figure}
Starting from the symplectic dual polar space $\ma{G}_n$, we construct its universal embedding 
$U(\ma{G}_n):=\FF_2(\ma{P}_n)/\eta(\FF_2(\ma{L}_n))$, where $\eta:\FF_2(\ma{L}_n)\longrightarrow\FF_2(\ma{P}_n)$ sends every line to the sum of its three elements.
Brouwer conjectured 
that the value of $\dim(U(\ma{G}_n))$ is given by the sequence $(x_n)_{n\in\NN} = (2,5,15,51,187,\ldots)$ with  $x_n=(2^n+1)(2^{n-1}+1)/3$ which is the sequence \href{https://oeis.org/A007581}{A007581} in \cite{oeis}. 
This conjecture was proved by P.~Li in \cite{li} and independently by A.~Blokhuis and A.~E.~Brouwer in \cite{brouwer}.
 In this paper we are mainly concerned with the procedure employed by P.~Li in \cite{li} where he considers sets $\ma{N}^n$ of subspaces of a fixed $n$-dimensional maximal isotropic subspace of a $2n$-dimensional symplectic space $V$ over $\FF_2$.
These subsets $\ma{N}^n$ are not maximal in a certain sense (see Section \ref{sec:embedding}). 
Every set $\ma{N}^n$ is subdivided into a disjoint union of families which are constructed inductively.
In our work we construct a bijection between the set $\ma{N}^n$ and a set of words of length $n$ in so-called standard order, formed by four different letters. 
Moreover, this bijection respects the inductive construction thus allowing us to construct every element of $\ma{N}^n$ in a very simple way. 
As a consequence, our procedure gives an alternative proof of the formula for the dimension of the universal embedding $U(\ma{G}_n)$ since the number of words can be easily counted.
This construction establishes a relationship between the first and the second of the many different interpretations of the sequence $(x_n)_{n\in\NN} = (2,5,15,51,187,\ldots)$ (the sequences \href{https://oeis.org/A007581}{A007581} and 
\href{https://oeis.org/A124303}{A124303}
in \cite{oeis}) in the following list:
 \begin{enumerate}
     \item The dimension of the universal embedding of the symplectic dual polar space \cite{brouwer, li}.
     \item The density of a language with four letters \cite{mor,SW2015}.
     \item The number of isomorphism classes of regular fourfold coverings  of a graph $L$ with Betti number $n = \beta(L)$ and with voltage group $\FF_2 \times \FF_2$ \cite{jin}.
     \item The number of non-equivalent states of a Hanoi graph associated to the Hanoi tower with $n$ discs and four pegs \cite{hinz}.
     \item The dimension of a certain centralizer algebra associated to a group of order 96 \cite{kosuda}.
     \item The dimension of the space of symmetric polynomials in 4 noncommuting variables \cite{rosas1,rosas2}.
     \item The cob-invariant associated to the group $\ZZ_2^n$ \cite{carlos1,carlos3}.
\end{enumerate}
Actually, all this is part of a more general setting with an arbitrary prime number $p$.
In \cite{SW2015} the authors consider a language with $p^2$ letters as a quotient of  $(\ZZ_p\times\ZZ_p)^n$ by the special linear group $SL(2,\ZZ)$. 
In the case of a dual polar space, we consider the totally isotropic subspaces of an $\FF_p$-vector space $V$, where we get configurations with points $\ma{P}_n$ and lines $\ma{L}_n$ satisfying

\begin{equation*}
   |\ma{P}_n(p)|=\prod_{k=1}^n(p^k+1)\,\textrm{ and }|\ma{L}_n(p)|=\frac{(p^n-1)|\ma{P}_n(p)|}{p^2-1}
\end{equation*}
where every line has $p+1$ points and every point is contained in $\frac{p^n-1}{p-1}$ lines, see Remark \ref{corq}. 
For instance, for $n=2$ this produces a sequence of self-dual configurations $\big( (p+1)(p^2+1) \big)_{p+1}$ for prime numbers $p$, i.e., $15_{2+1}$, $40_{3+1}$, $156_{5+1}$, $400_{7+1}, \dots$ which we will call the $p$-Cremona-Richmond configurations.
The sequence $15,40,156,400,1464,\dots$ appears as the sequence 
\href{https://oeis.org/A131991}{A131991} in \cite{oeis}.

The paper is organized as follows.
In Section~\ref{sec1} we discuss the set $W^n$ of words of length $n$ in so-called standard order formed by four letters and we give a procedure for the construction of all the words $W^{n+1}$ from the ones in $W^n$.
We obtain two proofs for the formula for $|W^n|$ (Proof of formula~\eqref{eq:formula} on page~\pageref{page:proof} and Remark~\ref{rem:cases}).
In Sections~\ref{sec:vectorspaces} and \ref{sec:symplectic} we prove several facts on isotropic subspaces of symplectic $\FF_2$-vector spaces and the symplectic dual polar space which seem to be well-known but the proofs are not easily accessible in the literature or scattered around.
In Section~\ref{sec:embedding} we review Li's proof for the formula of the dimension of the universal embedding of the symplectic polar space $\mathrm{Sp}_{2n}(2)$ which allows us to construct the bijection between the words $W^n$ and Li's vector spaces $\ma{N}^{n+1}$ in Section~\ref{sec:bijection}. 
This gives a new proof for the formula for the dimension of the universal embedding of the symplectic dual polar space in Theorem~\ref{thm:bijection}.
In Appendix~\ref{appendix:CR} we present the decomposition of the collinearity graph $\Gamma$ for $n=2$ and $n=3$ in its subgraphs $\Gamma_k$. 
In Appendix~\ref{appendix:construction} we show the construction of $W^{n+1}$ from $W^n$ for $n=1,2,3,4$.
Finally, in Appendix~\ref{apenC}, we present a classification of words in $W^n$ according to the 8 cases specified in Section \ref{sec1}.

\section{Languages}
\label{sec1}
Let us consider a language with the $4$ letters $0,1,2,3$.
For $n\in\NN$ we define $\widetilde W^n := \{ a_1\dots a_n : a_j=0,1,2,3\}$ to be the set of all possible words of length $n$ formed by the letters $0,1,2,3$.
In this article we will be mainly concerned with the subset $W^n$ of words in the so-called \define{standard order} \cite{ArndtSloane2016, mor}.
The set $W^n$ consists of the words $a_1 a_2 \dots a_n\in\widetilde W^n$ such that there exist $1\le j<k$ with:
\begin{itemize}
   \item[(R1)] $a_i = 0$\; for\; $i < j$,
   \item[(R2)] $a_j = 1$,
   \item[(R3)] $a_i \in\{0, 1\}$\; for\; $j<i<k$,
   \item[(R4)] $a_k = 2$\; if\; $k\le n$,
   \item[(R5)] $a_i \in\{0, 1, 2,3\}$\; for\; $i>k$.
\end{itemize}
Note that (R5) applies only if $k< n$.
For a word $a=a_1a_2\dots a_n$ the rules above can be written compactly as
\begin{equation*}
   \label{eq1} 
   0\le a_i\leq \op{max}_{j< i}\{a_j\}+1, 
   \qquad 1\le i \le n.
\end{equation*}
Note that the set $W^n$ is a special case for $p=2$, of 
the more general sets $W_p^n$ for $p=2$ defined in \cite{SW2015}.
\smallskip

\begin{definition} 
   The cardinality of $W^n$ is called the \define{density of the language} $W^n$.
   We use the notation $g_W(n) := |W^{n}|$.
\end{definition}
\begin{example} 
   \label{eq:words}
   We have
   $W^1=\{0\}$, consisting of $1$ word,
   $W^2=\{00, 01\}$, consisting of $2$ words,
   $W^3=\{000, 001, 010, 011, 012\}$, consisting of $5$ words.
   For $n=4$ there are 15 words and the elements of $W^4$ are 
\begin{center}   \begin{tabular}[t]{lllll}
      0000 & 0001 & 0010 & 0011 & 0012 \\
      0100 & 0101 & 0102 & 0110 & 0111 \\
      0112 & 0120 & 0121 & 0122 & 0123\,.
   \end{tabular}
\end{center}
\end{example}
The following theorem was already shown in \cite{mor} and \cite{SW2015}.
We will give an alternative proof after Proposition~\ref{prop:W}.
\begin{theorem}
   \label{thm:formula}
   For $n\in \NN_0$ the density of the language $W^{n}$ is
   \begin{equation}
      \label{eq:formula}
      |W^{n}| = g_W(n) = \frac{(2^{n-1}+1)(2^{n-2}+1)}{3}\,.
   \end{equation}
\end{theorem}

In the present work, we want to provide a different point of view, motivated by the work of Li \cite{li}, and we will give an alternative proof of Theorem~\ref{thm:formula}.
In what follows we study some facts which are fundamental for the proof of this theorem.
\bigskip

Let $n\ge 2$.
We will show how all words in $W^{n+1}$ can be constructed from the words in $W^{n}$.
\begin{itemize}
   \item{\bf Case 1.}
   Take an arbitrary word in $W^n$ and attach $0$ at the end. 
   This gives a valid word in $W^{n+1}$.
   The number of all such words is $g_W(n)$.

   \item{\bf Case 2.}
   Take an arbitrary word in $W^n$ and attach $1$ at the end. 
   This gives a valid word in $W^{n+1}$ and it is not contained in the words obtained in case 1.
   The number of all such words is $g_W(n)$.
   \end{itemize}

   \noindent
   For Cases 3, 4, 5 we take an arbitrary word $a=a_1a_2\dots a_n$ in $W^n$ which ends in $2$ or $3$.
   Note that this implies $1\in\{a_1,\, \dots,\, a_{n-1}\}$ and that
   therefore $a=a_1a_2\dots a_{n-1} \ell a_n$ is a valid word in $W^{n+1}$ for $\ell=0,1,2$.

   \begin{itemize}
   \item{\bf Case 3.}
   Insert the letter $0$ before $a_n$.
   We obtain the valid word $\widetilde a=a_1a_2\dots a_{n-1}0a_n\in W^{n+1}$.
   Clearly this word is not contained in the words constructed so far.

   \item{\bf Case 4.}
   Insert the letter $1$ before $a_n$.
   We obtain the valid word $\widetilde a=a_1a_2\dots a_{n-1}1a_n\in W^{n+1}$.
   Clearly this word is not contained in the words constructed so far.

   \item{\bf Case 5.}
   Insert the letter $2$ before $a_n$.
   We obtain the valid word $\widetilde a=a_1a_2\dots a_{n-1}2a_n\in W^{n+1}$.
   Clearly this word is not contained in the words constructed so far.
\end{itemize}

\noindent
The number of words in each of the Cases 3, 4, 5 is 
\begin{align*}
   \#(\text{words of length $n$ ending in } 2 \text{ or } 3)
   &= g_W(n)-\#(\text{words of length $n$ ending in } 0 \text{ or } 1)
   \\
   &= g_W(n) - 2 g_W(n-1),
\end{align*}
because
$ \#(\text{words of length $n$ ending in } 0 )
= \#(\text{words of length $n$ ending in } 1)
= g_W(n-1)$ as in Case~1 and Case~2.

\begin{itemize}
   \item{\bf Case 6.}
   Let $a=a_1a_2\dots a_n$ in $W^n$ which ends in $2$ or $3$
   and such that $2\in\{a_1,\, \dots,\, a_{n-1}\}$.\\
   Insert the letter $3$ before $a_n$.
   We obtain the valid word $\widetilde a=a_1a_2\dots a_{n-1}3a_n\in W^{n+1}$.
   Clearly this word is not contained in the words constructed so far.

   \item{\bf Case 7.}
   Let $a=a_1a_2\dots a_n$ in $W^n$ which ends in $2$ or $3$
   and such that $2\not\in\{a_1,\, \dots,\, a_{n-1}\}$.\\
   This implies that $a_n=2$ and $a_j\in\{0,1\}$ for $1\le j \le n-1$.
   ttach $3$ to obtain the new word
   $\widetilde a=a_1a_2\dots a_{n-1}23\in W^{n+1}$.
   Clearly this word is not contained in the words constructed so far.
   The number of all such words is equal to the number of strings of length $n-1$ consisting only of $0$ and $1$, with exception of the zero string.
   So the number of the words in this case is $2^{n-2}-1$.
\end{itemize}

\noindent
   The total number of words in the Cases~6 and 7 together is
   \begin{align*}
   \#(\text{words of length $n$ ending in } 2 \text{ or } 3)
   &= g_W(n)-\#(\text{words of length $n$ ending in } 0 \text{ or } 1)
   \\
   &= g_W(n) - 2 g_W(n-1).
   \end{align*}

\begin{itemize}
   \item{\bf Case 8.}
   Let 
   $\widetilde a=0\dots 0 1 2\in W^{n+1}$.
   Clearly this word is not contained in the words constructed so far.
\end{itemize}

We say that a word $\widetilde a\in W^{n+1}$ \define{is in Case~k} for $k=1,\,\dots,\, 8$, if it is constructed from a word $a\in W^{n}$ as described in Case~$k$.
\bigskip

\begin{prop}
   \label{prop:W} 
   Let $n\ge 2$. 
   Then each word in $W^{n+1}$ is constructed as in exactly one of the Cases 1 -- 8 above.
\end{prop}
\begin{proof}
   Let $\widetilde a=a_1a_2\dots a_{n}a_{n+1}\in W^{n+1}$.

   \begin{itemize}
      \item
      Suppose that $a_{n+1}\in \{0,1\}$. 
      Then clearly $\widetilde a$ is constructed either as in Case~1 or in Case~2.

      \item
      Suppose that $a_{n+1}\in \{2,3\}$. 
      Note that this implies $1\in\{a_1, \dots, a_{n}\}$. 
      \begin{itemize}
	 \item If $a_n=0$, we can erase it and obtain the valid word 
	 $a=a_1\dots a_{n-1}a_{n+1}\in W^n$.
	 If we now apply the procedure of Case 3, we recover $\widetilde a$.

	 \item If $a_n=1$ and $a=a_1\dots a_{n-1}a_{n+1}$ is a valid word in $W^n$,
	 then we can apply the procedure of Case 4 and we obtain again $\widetilde a$.
	 If $a=a_1\dots a_{n-1}a_{n+1}$ is not a valid word in $W^n$, then necessarily
	 $\widetilde a= 0\dots 012$ and we have the word of Case~8.

	 \item If $a_n=2$ and $a=a_1\dots a_{n-1}a_{n+1}$ is a valid word in $W^n$,
	 then we can apply the procedure of Case 5 and we obtain again $\widetilde a$.
	 If $a=a_1\dots a_{n-1}a_{n+1}$ is not a valid word in $W^n$, then necessarily
	 $a_{n+1}=3$ and $2\not\in\{a_1,\,\dots,\, a_{n-1}\}$.
	 Then we can apply the procedure of Case 7 to the word
	 $a'=a_1\dots a_{n-1}a_{n}$ and we recover $\widetilde a$.

	 \item If $a_n=3$ then necessarily $a=a_1\dots a_{n-1}a_{n+1}$ is a valid word in $W^n$
	 and we can apply the procedure of Case 6 to recover $\widetilde a$.
	 \qedhere

      \end{itemize}
   \end{itemize}
\end{proof}

In Appendix~\ref{apenC} we give a classification of the words of $W^n$ according to the cases described above and 
in Appendix~\ref{appendix:construction} 
we show explicitly how $W^{n+1}$ is constructed from $W^n$ for $n=2,3,4$.
Proposition~\ref{prop:W} allows us to prove the formula \eqref{eq:formula} as follows.

\begin{proof}[Proof of Formula~\eqref{eq:formula}]
   \label{page:proof} 
   By Proposition~\ref{prop:W} we know that, for $n\ge 2$,
   \begin{align*}
      |W^{n+1}| = g_W(n+1)= 2g_W(n) + 4[ g_W(n)-2g_W(n-1)] + 1 = 6 g_W(n) - 8 g_W(n-1) + 1.
   \end{align*}

   From Example~\ref{eq:words} we obtain that 
   $|W^1|=1$, $|W^2|=2$, 
   hence formula~\eqref{eq:formula} is satisfied for $n=1,2$.
   Now suppose that the formula holds for all $j\le n$. Then
   \begin{align*}
      g_W(n+1) & = 6 g_W(n) - 8 g_W(n-1) + 1
      =6\frac{  (2^{n-1}+1)  (2^{n-2}+1) }{3} - 8\frac{  (2^{n-2}+1) (2^{n-3}+1) }{3}+1\\
      &=\frac{ (2^{n-2}+1)  }{3} \left[   6(2^{n-1}+1)  -8(2^{n-3}+1)   \right]+1
      \\
      &=\frac{ 2(2^{n-2}+1)(2^n-1) + 3 }{3} \\
      &=\frac{ 2^{2n-1} + 2^{n+1} - 2^{n-1} + 1}{3}
      =\frac{ 2^{2n-1} + 2^{n} + 2^{n-1} + 1}{3}\\
      &=\frac{(2^{n}+1)(2^{n-1}+1)}{3}.
      \qedhere
   \end{align*}
\end{proof}

\begin{remark}
   An alternative proof of Proposition~\ref{prop:W} makes use of the formula~\eqref{eq:formula} for $|g_W(n)|$ which was already proved in \cite{mor} and \cite{SW2015}.
   Then it is sufficient to prove that the number of words obtained by the Cases 1 to 8 is equal to $|W^{n+1}|$ because we already know that all cases are disjoint and that every word constructed in these cases belongs to $W^{n+1}$.
   That is, we have to show that
   \begin{equation*}
      g_W(n+1) = 2g_W(n) + 6[ g_W(n)-g_W(n-1) ] + 1.
   \end{equation*}
   This is a straightforward calculation.
\end{remark}
\smallskip

Proposition~\ref{prop:W} gives yet another way to calculate $|W^n|$ as the next remark shows.

\begin{remark}
   \label{rem:cases}
   We showed that with the rules in Cases~1 to 8, each word in $W^n$ which ends in $0$ or $1$ gives rise to exactly two words in $W^{n+1}$ (Cases 1 and 2). 
   They again end in $0$ or $1$.
   Each word in $W^n$ which ends in $2$ or $3$ gives rise to exactly six words in $W^{n+1}$, two of which end in $0$ or $1$, and four of them end again in $2$ or $3$ (Cases 3, 4, 5 and either 6 or 7).
   In addition we have the word from Case 8.
   \smallskip

   \noindent
   Let $g_W(n)=|W^n|$ and
   \begin{align*}
      s(n) &= \text{number of words in } W^n \text{ which end in } 0 \text{ or } 1,\\
      t(n) &= \text{number of words in } W^n \text{ which end in } 2 \text{ or } 3.
   \end{align*}
   Then we obtain
   $t(1)=t(2)=0,\ t(3)=1$, 
   $s(1)=1,\ s(2)=2,\ s(3)=4$ and
   \begin{align*}
      g_W(n+1) = 2 s(n) + 6 t(n) + 1,\quad
      s(n+1) = 2 s(n) + 2 t(n),\quad
      t(n+1) = 4 t(n) + 1,
      \qquad n\ge 2.
   \end{align*}
   Iterating the formula for $t(n)$, we find $t(n)=\frac{4^{n-2}-1}{3}$.
   For $s(n)$ we find
   \begin{align*}
      s(n) &= 2^{n-3}s(3) + \sum_{j=1}^{n-3} 2^j t(n-j)
      = 2^{n-1} + \frac{1}{3} \sum_{j=1}^{n-3} 2^j (2^{2n-4-2j}-1)
      = 2^{n-1} + \frac{2}{3} \left(\sum_{j=0}^{n-4} 2^{2n-6-j} -\sum_{j=0}^{n-4} 2^j \right)
      \\
      &= 2^{n-1} + \frac{2}{3} \left(2^{n-2}-1 \right) \sum_{j=0}^{n-4} 2^j
      = 2^{n-1} + \frac{2}{3} ( 2^{n-2}-1 )( 2^{n-3}-1 )
      \\
      &= \frac{1}{3} ( 3\cdot 2^{n-1} + 2^{2n-4} - 2^{n-1} - 2^{n-2} + 2 )
      = \frac{1}{3} ( 2^{2n-4} + 2^{n-1} + 2^{n-2} + 2 ).
   \end{align*}
   So we find again formula~\eqref{eq:formula} for $g_W(n)$:
   \begin{align*}
      g_W(n) &= s(n) + t(n) 
      = \frac{1}{3} \left[ 
      2^{2n-4} + 2^{n-1} + 2^{n-2} + 2
      + 4^{n-2} - 1 \right]
      = \frac{1}{3} \left[ 
      2^{2n-3} + 2^{n-1} + 2^{n-2} + 1
      \right]
      \\
      &= \frac{(2^{n-1}+1)(2^{n-2}+1)}{3}.
   \end{align*}
\end{remark}
\smallskip

\begin{definition}
Let us introduce some more notation.
We define the following subsets of $W^n$:
\begin{align*}
   W_0^n &:= \{ a_1\dots a_{n-1}0: a_j=0,1,2,3\} = \text{ all words ending in } 0,\\
   W_1^n &:= \{ a_1\dots a_{n-1}1: a_j=0,1,2,3\} = \text{ all words ending in } 1,\\
   S_0^n &:= \{ a_1\dots a_{n-2}0b_n: a_j=0,1,2,3,\ b_n=2,3\} = \text{ all words ending in $2$ or $3$ with } a_{n-1}=0,\\
   S_1^n &:= \{ a_1\dots a_{n-2}1b_n: a_j=0,1,2,3,\ b_n=2,3\} = \text{ all words ending in $2$ or $3$ with } a_{n-1}=1,\\
   S_2^n &:= \{ a_1\dots a_{n-2}2b_n: a_j=0,1,2,3,\ b_n=2,3\} = \text{ all words ending in $2$ or $3$ with } a_{n-1}=2,\\
   S_3^n &:= \{ a_1\dots a_{n-2}3b_n: a_j=0,1,2,3,\ b_n=2,3\} = \text{ all words ending in $2$ or $3$ with } a_{n-1}=3,\\
   S^n &:= S_0^n \cup S_1^n \cup S_2^n \cup S_3^n\\
   C_1^n &:= \{ a_1\dots a_{n-2}23: a_j=0,1\},\\
   C_2^n &:= \{ a_1\dots a_{n-1}2: a_j=0,1\} = \text{ all words with $a_n=2$ and no other $2$ },\\
   C_3^n &:= \{ 0\dots 012\}.
\end{align*}
\noindent
Observe that $C_1^n\subset S_2^n$,
$C_3^n\subset S_1^n$
and that 
$W^{n+1} =
W_0^{n+1} \sqcup W_1^{n+1} \sqcup
S_0^{n+1} \sqcup S_1^{n+1} \sqcup
S_2^{n+1} \sqcup S_3^{n+1}$,
where $\sqcup$ denotes the disjoint union.
Now we define the \define{insert operators} for $k=1,\dots, n+1$ and $\ell=0,1,2,3$ as follows:
\begin{equation*}
   A_{k,\ell}^n: W^n\to\widetilde W^{n+1},
   \qquad
   A_{k,\ell}^n(a_1\dots a_n) = a_1\dots a_{k-1} \ell a_{k}\dots a_n.
\end{equation*}
For $n\in\NN$ and $1\le j \le n$ we define the \define{erase operators}
\begin{equation*}
   E_{j}^n: W^{n}\to \widetilde W^{n-1},
   \qquad
   E_{j}^n(a_1\dots a_n) = 
   a_1\dots a_{j-1}a_{j+1}\dots a_{n}.
\end{equation*}
\end{definition}
It should be observed that for $a\in W^n$ and $j\in\{1,2,\dots,n\}$, the word $E_j^n(a)$ is not necessarily a word in $W^{n-1}$. 
\bigskip

\noindent
With this new notation, the results of this section so far can be summarized as follows.

\begin{theorem}
   \label{thm:summary1}
   Let $n\in\NN$.
   Then the following maps are bijections:

   \begin{minipage}{.35\textwidth}
      \begin{alignat*}{2}
	 \tag{Case 1}
	 A^n_{n+1,0} &:\ &W^n&\to W^{n+1}_0,\\
	 \tag{Case 2}
	 A^n_{n+1,1} &:\ &W^n&\to W^{n+1}_1,\\
	 \tag{Case 3}
	 A^n_{n,0} &:\ &S^n&\to S^{n+1}_0,
      \end{alignat*}
   \end{minipage}
   \hspace{\fill}
   \begin{minipage}{.5\textwidth}
      \begin{alignat*}{2}
      \tag{Case 4}
      A^n_{n,1} &:\ &S^n&\to S^{n+1}_1\setminus C^{n+1}_3,\\
      \tag{Case 5}
      A^n_{n,2} &:\ &S^n&\to S^{n+1}_2\setminus C^{n+1}_1,\\
      \tag{Case 6}
      A^n_{n,3} &:\ &S^n\setminus C^{n}_2&\to S^{n+1}_3,\\
      \tag{Case 7}
      A^n_{n+1,3} &:\ &C^{n}_2&\to C^{n+1}_1
   \end{alignat*}
   \end{minipage}
   \bigskip

   \noindent
   and
   $W^{n+1} =
   W_0^{n+1} \sqcup W_1^{n+1} \sqcup
   S_0^{n+1} \sqcup \big(S_1^{n+1}\setminus C_3^{n+1}\big) \sqcup
   \big(S_2^{n+1}\setminus C_1^{n+1}\big) \sqcup
   S_3^{n+1}\sqcup C_1^{n+1}\sqcup C_3^{n+1}$
   is the disjoint union of the images of the maps above and $C_3^{n+1}$.

   \noindent
   The inverses of the maps above are 

   \begin{minipage}{.35\textwidth}
      \begin{alignat*}{2}
	 E^{n+1}_{n+1} &:\ & W^{n+1}_0 &\to W^n,\\
	 E^{n+1}_{n+1} &:\ & W^{n+1}_1 &\to W^n,\\
	 E^{n+1}_{n} &:\ & S^{n+1}_0 &\to S^{n},
      \end{alignat*}
   \end{minipage}
   \hspace{\fill}
   \begin{minipage}{.5\textwidth}
      \begin{alignat*}{2}
	 E^{n+1}_{n} &:\ &S^{n+1}_1\setminus C^{n+1}_3&\to S^n,\\
	 E^{n+1}_{n} &:\ &S^{n+1}_2\setminus C^{n+1}_1&\to S^n,\\
	 E^{n+1}_{n} &:\ &S^{n+1}_3&\to S^n\setminus C^{n}_2,\\
	 E^{n+1}_{n+1} &:\ &C^{n+1}_1&\to C^{n}_2.
   \end{alignat*}
   \end{minipage}

   \noindent
   Moreover,
   \begin{align*}
   |W^{n+1}_0| = |W^{n+1}_1| &= g_W(n), \\
   |S^{n+1}_0| = |S^{n+1}_1\setminus C^{n+1}_3| = |S^{n+1}_2\setminus C^{n+1}_1| &= g_W(n)-2g_W(n-1),\\
   |S^{n+1}_3| + |C^{n+1}_1| &= g_W(n)-2g_W(n-1).
   \end{align*}

\end{theorem}

\begin{theorem}
   \label{thm:wordsconstuction}
   For every word $a=a_1\dots a_n\in W^n$ exactly one of the following holds.
   \begin{enumerate}
      \item 
      There is exactly one sequence of maps $A^1, A^2, \dots, A^{n-1}$ such that 
      $a=A^{n-1} \cdots A^{1}(0)$
      where the $A^j$ are maps of type $A^j_{\ell, a}$ as in Theorem~\ref{thm:summary1}.

      \item 
      There is exactly one $k\le n$ and exactly one sequence of maps $A^k, A^{k+1}, \dots, A^{n-1}$ such that 
      $a=A^{n-1} \cdots A^{k}(a')$
      where the $A^j$ are maps of type $A^j_{\ell, a}$ as in Theorem~\ref{thm:summary1} and $a'=0\dots 012\in W^k$.
   \end{enumerate}
\end{theorem}


\section{Isotropic subspaces} 
\label{sec:vectorspaces}

Let $n$ be an integer and consider the vector space $V=\FF_2^{2n}$ over the field with two elements $\FF_2$ equipped with a symplectic form $\iprod{\cdot}{\cdot}: V\times V\to \FF_2$, that is, a non-degenerate bilinear form which satisfies
\begin{equation}
   \label{eq:sym}
   \iprod{u}{v}= \iprod{v}{u}
   \quad\text{and}\quad \iprod{v}{v}=0
   \qquad \text{for all }\ u,v\in V.
\end{equation}
Non-degeneracy means that for every $v\in V\setminus\{0\}$ there exists $w\in V$ such that $\iprod{v}{w}\neq 0$.
A basis $u_1,\dots, u_n, v_1,\dots, v_n$ of $V$ is called a \define{symplectic basis of $V$} if 
   \begin{equation}
      \label{eq:symplecticbase}
      \iprod{u_j}{u_k} = \iprod{v_j}{v_k} =0,
      \qquad
      \iprod{u_j}{v_k} = \delta_{jk},
      \qquad j,k=1,\dots, n.
   \end{equation}
For a subspace $U$ of $V$ we define its \define{orthogonal complement} by 
\begin{equation*}
   U^\perp := \{ x\in V: \iprod{x}{u}=0 \text{ for all } u\in U\}.
\end{equation*}
We write $U\perp W$ if $[u,w]=0$ for all $u\in U$ and $w\in W$.
Moreover, we write $U\operp W := U\oplus W$ if $U$ and $W$ are subspaces with $U\cap W = \{0\}$ and $U\perp W$.

Observe that the non-degeneracy of $\iprod{\cdot}{\cdot}$ implies that the natural map $\Phi:V\to V^*,\ v\mapsto \iprod{\cdot}{v}$ from $V$ into its dual space $V^*$ is a bijection.
In particular, it maps $U^\perp$ bijectively to the $(\dim V - \dim U)$-dimensional subspace $\{\varphi\in V^* : \varphi(u)=0 \text{ for all } u\in U\}$.
It follows that 
\begin{equation}
   \label{eq:ComplSubspace}
   \dim U + \dim U^\perp = \dim V = 2n.
\end{equation}
A subspace $U$ of $V$ is called \define{totally isotropic} if $U\subseteq U^\perp$, that is, if $\iprod{v}{w}=0$ for every pair of elements $v,w\in U$. 
A subspace $U$ of $V$ is called a \define{maximal totally isotropic subspace} if it is totally isotropic and not properly contained in any other totally isotropic subspace of $V$.
This is the case if and only if $U=U^\perp$.
Every isotropic subspace is contained in a maximal totally isotropic subspace.
Since every one-dimensional subspace of $V$ is isotropic by~\eqref{eq:sym}, it follows that maximal totally isotropic subspaces exist.
Although the following results on (totally) isotropic subspaces seem to be well-known, we include their proofs here.
\begin{prop}
   \label{prop:max}
   Every maximal totally isotropic subspace $U$ of $V$ has dimension $n$.
   Moreover, every basis $u_1,\dots, u_n$ of $U$ can be completed to a symplectic basis of $V$.
\end{prop}

\begin{proof}
   If $U$ is maximal totally isotropic, then $U = U^\perp$ and it follows from~\eqref{eq:ComplSubspace}
   that $\dim U = \frac{1}{2} \dim V=n$.
   Let $u_1,\dots, u_n$ be a basis of the maximal totally isotropic subspace $U$
   and set $U_1:= \linspan\{u_2,\dots, u_n\}$.
   Then clearly $\dim U_1=n-1$, 
   $\dim U_1^\perp=n+1$ and 
   $U\subseteq U_1^\perp$.
   Assume that $\iprod{x}{u_1}=0$ for all $x\in U_1^\perp$.
   Then $U_1^\perp \subseteq U^\perp$ which contradicts
   $\dim U_1^\perp = n+1 > n = \dim U^\perp$.
   Therefore, we can choose $w_1\in U_1^\perp$ such that $\iprod{u_j}{w_1}=\delta_{j1}$.
   Note that $\iprod{w_1}{w_1}=0$ by~\eqref{eq:sym}.

   Now set $V_1=\linspan\{u_1,w_1\}$ and $W_1=\linspan\{u_1,w_1\}^\perp$.
   Clearly, $\dim V_1 = 2$ and $V_1\cap W_1 = \{0\}$, hence
   \begin{equation}\label{eq:V1}
      V= V_1\operp W_1.
   \end{equation}
   If we restrict $\iprod{\cdot}{\cdot}$ to $W_1$, then $W_1$ becomes a nondegenerate symplectic space of dimension $2n-2$.
   Note that $U_1\subseteq W_1$ is a maximal totally isotropic subspace of $W_1$.

   Repeating the construction above $n-1$ times, we find vectors $w_1,\dots, w_n$ with $[u_j,w_k]=\delta_{jk}$ and 
   \begin{equation*}
      V=V_1 \oplus \cdots\oplus V_n =V_1 \operp \cdots\operp V_n
   \end{equation*}
   for $V_j := \linspan\{ u_j,w_j\}$, $j=1,\dots, n$.
\end{proof}

\begin{remark}
   The proposition shows that, after choosing a symplectic basis of $V$,
   the symplectic form $\iprod{\cdot}{\cdot}$ can be written as
   \begin{equation*}
      \iprod{x}{y}  := xJy^t := \scalar{x}{Jy} = \scalar{Jx}{y}
      \quad\text{with}\quad
      J=
      \begin{pmatrix}
	 0_n & I_n \\
	 I_n & 0_n \\
      \end{pmatrix}
   \end{equation*}
   where $0_n$ and $I_n$ are the $n\times n$ zero matrix and the $n\times n$ identity matrix, respectively, $x,y$ are given as column vectors in the symplectic basis
   and $\scalar{\cdot}{\cdot}$ is the ``euclidian inner product'' with respect to the basis
   (that is, $\scalar{u_j}{v_k} = 0$, $\scalar{u_j}{u_k} = \scalar{v_j}{v_k} = \delta_{jk}$.)
\end{remark}

\begin{theorem}
   \label{thm:SympBasis}
   Let $V$ be a symplectic vector space of dimension $2n$ and let $X,P\subset V$ be maximal totally isotropic subspaces of $V$.
   Then $\dim(X\cap P)=n-k$ if and only if there exists a symplectic basis 
   $x_1,\dots, x_n, y_1,\dots, y_n$ of $V$ and numbers
   $a_{jm}\in\{0,1\}$ with $a_{jm}=a_{mj}$ for $j, m =1,\dots, k$ such that
   \begin{align}
      \label{eq:xp}
      X=\linspan\{x_1,\dots, x_n\},
      \qquad
      P=\linspan\{p_1,\dots,p_k, x_{k+1},\dots, x_n\}
      \quad \text{ with }\ p_j= y_j + \sum_{\ell=1}^k a_{j\ell}x_\ell.
   \end{align}
\end{theorem}
\begin{proof}
   Suppose that $P$ is as in \eqref{eq:xp}.
   Then it is easy to see that $P$ is an isotropic subspace of dimension $n$ and that 
   $X\cap P= \linspan\{x_{k+1},\dots, x_n\}$.
   \smallskip

   Now assume that $X$ and $P$ are maximal totally isotropic subspaces of $V$ with $\dim(X\cap P)=n-k$ for some $k=0,\dots, n$.
   Let us choose a basis $x_{1}, \dots, x_n$
   of $X$ such that 
   $X\cap P = \linspan\{ x_{k+1},\dots, x_{n}\}$ 
   and complete it to a symplectic basis
   $x_1,\dots, x_n, y_1,\dots, y_n$ of $V$.
   Now we choose $z_1,\dots, z_k\in P$ such that 
   $P = \linspan\{ z_1,\dots, z_{k}, x_{k+1},\dots, x_n\}$.
   Clearly they can be chosen such that
   \begin{equation*}
      z_j =  \sum_{\ell=1}^n b_{j\ell}y_\ell + \sum_{\ell=1}^k a_{j\ell}x_\ell
   \end{equation*}
   for $a_{j\ell},\, b_{j\ell}\in\{0,1\}$.
   Since $P$ is totally isotropic, we must have for all $j=1,\dots, k$ and $m=k+1,\dots, n$
   \begin{equation*}
      0 = [z_j,\, x_m]
      = \left[ 
      \sum_{\ell=1}^n b_{j\ell}y_\ell
      + \sum_{\ell=1}^k a_{j\ell}x_\ell,\ x_m 
      \right]
      = \sum_{\ell=1}^n b_{j\ell} [y_\ell,\ x_m ]
      + \sum_{\ell=1}^k a_{j\ell} [ x_\ell, x_m] 
      = b_{jm},
   \end{equation*}
   hence $z_j = \sum_{\ell=1}^k b_{j\ell}y_\ell + \sum_{\ell=1}^k a_{j\ell}x_\ell$. 
   By taking appropriate linear combinations if necessary, we may assume without restriction that 
   $z_j = b_{jj}y_j + \sum_{\ell=1}^k a_{j\ell}x_\ell$.
   If one of the $b_{jj}$ were $0$, then the corresponding $z_j$ would belong to $X$ and therefore we would have $\dim(X\cap P)> n-k$.
   So we may assume without restriction that
   \begin{equation*}
      z_j = y_j + \sum_{\ell=1}^k a_{j\ell}x_\ell,
      \qquad
      j=1,\dots, k.
   \end{equation*}
   Using again that $P$ is totally isotropic, we obtain for $j, m\in \{1,\dots, k\}$ that
   \begin{equation*}
      0 = [z_j,\, z_m]
      = \left[ y_j + \sum_{\ell=1}^k a_{j\ell}x_\ell,\ y_m  + \sum_{\ell=1}^k a_{m\ell}x_\ell \right]
      = \left[ \sum_{\ell=1}^k a_{j\ell}x_\ell,\ y_m  \right]
      + \left[ y_j,\ \sum_{\ell=1}^k a_{m\ell}x_\ell \right]
      = a_{jm} + a_{mj},
   \end{equation*}
   which implies $a_{jm} = a_{mj}$ for $j, m\in \{1,\dots, k\}$.
\end{proof}

Note that the representation of $P$ as in \eqref{eq:xp} is unique.

\begin{cor}
   \label{cor:N}
   Let $X$ be a maximal totally isotropic subspace of $V$ and let $\widetilde X\subseteq X$ be a subspace with $\dim \widetilde X = n-k$ for some $k\in\{0,\dots, n\}$.
   Then there exist exactly $2^{\frac{1}{2}k(k+1)}$ maximal totally isotropic subspaces $P$ with $X\cap P = \widetilde X$.
\end{cor}
\begin{proof}
   Since every maximal totally isotropic subspace $P$ with $X\cap P=\widetilde X$ must be of the form \eqref{eq:xp}, the number of such subspaces is equal to the number of all possible choices of $a_{jm}\in\{0,1\}$ with $1\le j \le m\le k$ (which is equal to the number of all possible symmetric $(k\times k)$-matrices over $\FF_2$).
   This number is 
   $2^{\sum_{\ell=1}^k \ell} = 2^{\frac{1}{2}k(k+1)}$.
\end{proof}

In the following we use the \define{$p$-binomial coefficients} defined by
\begin{equation*}
   \binom{n}{k}_p = \frac{ [n]_p !}{ [k]_p! [n-k]_p !}
   \quad\text{where }\
   [n]_p!=[n]_p[n-1]_p\cdots [1]_p
   \ \text{ and }\
   [n]_p = 1 + p + p^2 + \cdots + p^{n-1} = \frac{p^n-1}{p-1}.
\end{equation*}
In this case, $p$ is $2$.

\begin{cor}\label{cor:contained}
   Every totally isotropic subspace $\widetilde X$ with $\dim\widetilde X = n-k$ is contained in exactly 
   \begin{equation*}
      \sum_{j=0}^k 
      \binom{k}{j}_2
      2^{\frac{1}{2} j(j+1)}\,
   \end{equation*}
   different maximal totally isotropic subspaces.
\end{cor}
\begin{proof}
   Let $X$ be a maximal totally isotropic subspace with $\widetilde X \subset X$.
   We obtain all maximal totally isotropic subspaces $P$ which contain $\widetilde X$, by first extending $\widetilde X$ to a subspace $Y$ with $\widetilde X \subseteq Y \subseteq X$ and then extending $Y$ to a maximal totally isotropic subspace $P$ with $P\cap X = Y$.
   If $\dim Y = n-k+\ell$, then, by Corollary~\ref{cor:N}, there are exactly
   $2^{\frac{1}{2} (k-\ell)(k-\ell+1)}$ such extensions.
   Therefore the number of all maximal isotropic extensions of $\widetilde X$ is
   \begin{equation}
      \label{eq:Nextensions}
      N :=
      \sum_{\ell=0}^k 
      \#\{\widetilde X \subseteq Y \subseteq X : \dim Y = n-k+\ell\}
      \cdot
      2^{\frac{1}{2} (k-\ell)(k-\ell+1)}.
   \end{equation}
   Clearly there exists a bijection between 
   the $(n-k+\ell)$-dimensional subspaces $Y$ of $X$ which contain $\widetilde X$
   and
   the $\ell$-dimensional subspaces of $X/\widetilde X$.
   This number is equal to $\binom{k}{\ell}_2$, see \cite[p.~28]{stanley}.
   So we obtain from \eqref{eq:Nextensions}
   \begin{equation*}
      N =
      \sum_{\ell=0}^k 
      \binom{k}{\ell}_2
      2^{\frac{1}{2} (k-\ell)(k-\ell+1)}
      =
      \sum_{j=0}^k 
      \binom{k}{k-j}_2
      2^{\frac{1}{2} j(j+1)}
      =
      \sum_{j=0}^k 
      \binom{k}{j}_2
      2^{\frac{1}{2} j(j+1)}.
      \qedhere
   \end{equation*}
\end{proof}

For the special case $k=1$ we obtain:
\begin{cor}
   \label{cor:tres}
   Every totally isotropic subspace $\widetilde X$ with $\dim\widetilde X = n-1$ is contained in exactly 3 different maximal totally isotropic subspaces.
   If $x_1,\dots, x_n,y_1,\dots, y_n$ is a symplectic basis of $V$ such that $\widetilde X = \linspan\{x_1,\dots, x_{n-1}\}$, then the maximal totally isotropic subspaces containing $\widetilde X$ are
   \begin{equation}
      \label{eq:3Lagrangians}
      \widetilde X \oplus
      \linspan\{x_n\},
      \qquad
      \widetilde X \oplus
      \linspan\{y_n\},
      \qquad
      \widetilde X \oplus
      \linspan\{x_n+y_n\}.
   \end{equation}
\end{cor}

In \cite[Lemma 9.4.1]{BCN},  Brouwer, Cohen and Neumaier give a formula for the number 
of all maximal totally isotropic subspaces of a $2n$-dimensional symplectic space.
In the following corollary we give an alternative proof for their formula.

\begin{cor}[Number of maximal totally isotropic subspaces]
   \label{cor:NN}
   The number of different maximal totally isotropic subspaces of a $2n$-dimensional symplectic space $V$ is
   \begin{equation*}
      \#(\text{maximal totally isotropic subspaces of } V) 
      = \prod_{k=1}^n (2^k+1).
   \end{equation*}
\end{cor}
\begin{proof}
   Let $X_0$ be a maximal totally isotropic subspace of $V$.
   It is known
   \cite{stanley}
   that the number of the $k$-dimensional subspaces of $X_0$ is $\binom{n}{k}_2$.
   Hence, by Corollary~\ref{cor:N}, the number of the maximal isotropic subspaces of $V$ with $k$-dimensional intersection with $X_0$ is 
   $\binom{n}{k}_2 2^{\frac{1}{2} k(k+1)}$. 
   Summation over $k=0,\dots, n$ gives 
   \begin{equation*}
      \#(\text{maximal totally isotropic subspaces of } V) 
      = \sum_{k=0}^n 
      \left[ \binom{n}{k}_2\ 2^{ \frac{1}{2} k(k+1) } \right].
   \end{equation*}
   Recall Gauss' binomial formula, see e.g. \cite[Ch. 5]{Kac},
   \begin{equation}
      \label{wer1}
      \prod_{j=1}^n(1+p^{j-1}x)=\sum_{k=0}^n\binom{n}{k}_pp^{\frac{1}{2}k(k-1)}x^k\, .
   \end{equation}
   Substituting $x=2$ and $p=2$ in \eqref{wer1} we get the result.
\end{proof}

In Corollary~\ref{cor:NHP} we will calculate the number of isotropic subspaces of dimension $n-1$. 
\medskip

For an arbitrary prime number $p$ we can show the following result.

\begin{remark}
\label{corq}For the field $\FF_p$, with $p$ a prime number, let $V$ be a $2n$-dimensional symplectic vector space over $\FF_p$. We have the following:
\begin{itemize}
    \item every totally isotropic subspace $\widetilde X$ with $\dim\widetilde X = n-k$ is contained in exactly 
   $\sum_{j=0}^k 
   \binom{k}{j}_p
   p^{\frac{1}{2} j(j+1)}$ different maximal totally isotropic subspaces; 
   \item the number of different maximal totally isotropic subspaces is 
\begin{equation*}
\sum_{k=0}^n 
      \left[ \binom{n}{k}_p\ p^{ \frac{1}{2} k(k+1) } \right]=\prod_{k=1}^n (p^k+1)\,;
\end{equation*}
\item the number of totally isotropic $(n-1)$-dimensional subspaces of $V$ is 
\begin{equation*}
   \frac{(p^n-1)}{p^2-1}\prod_{k=1}^n (p^k+1)
\end{equation*}
\end{itemize}

\end{remark} 

\section{Symplectic dual polar spaces}
\label{sec:symplectic}

A \emph{partial linear space} is a space $\mathcal G = (\mathcal L, \mathcal P)$ consisting of \define{points} $P\in\mathcal P$ and \define{lines} $L\in\mathcal L$ such that any line is incident with at least two points, and any pair of distinct points is incident with at most one line.
Two points $P,Q\in\mathcal P$ are called \emph{collinear} if they are on the same line $L\in\mathcal L$.
\smallskip

Let $V$ be a symplectic space of dimension $2n$.
It is easy to see that we obtain a partial linear space 
$\mG_n:=(\mP_n,\mL_n)$, called the \define{symplectic dual polar space}, if we set
   \begin{itemize}
      \addtolength{\itemsep}{-.4\baselineskip}
      \item $\mP_n = $ the set of all maximal totally isotropic subspaces of $V$;

      \item $\mL_n = $ the set of all totally isotropic subspaces of dimension $n-1$ of $V$.
   \end{itemize}
   Observe that in $\mathcal G_n$ every line is incident with exactly three points by Corollary~\ref{cor:tres}.
\smallskip

\begin{example}
   \label{ex:VS}
   Following Li, we write vectors $w\in V$ as row vectors $w=(w_1, \dots, w_{2n})$.
   For vectors $v_1,\dots, v_k\in V$, 
   we set
   $\left(
   \begin{smallmatrix} v_1\\ \vdots\\ v_k
   \end{smallmatrix}
   \right)
   := \linspan\{ v_1, \dots, v_k\}
   $.

   \begin{itemize}
      \item
      Let $n=1$. Then $V=\FF_2^2$ and its maximal isotropic subspaces are exactly the spans of the non-zero vectors of $V$.
      So we have
      $\mP_1=\{ (10),(01),(11)\}$ and $\mL_1 = \{(00)\}$, in particular $\op{dim}U(\mG_1)=2$.

      \item
      Let $n=2$. Then $V=\FF_2^4$ and the elements of $\mP_2$ are exactly the 15 following two-dimensional maximal isotropic subspaces:
      
      \begin{alignat*}{6} 
	 &\begin{pmatrix} 
	 1 0 0 0\\ 0 1 0 0
	 \end{pmatrix},\quad
	 &&\begin{pmatrix} 
	 1 0 0 0\\ 0 1 0 1
	 \end{pmatrix},\quad
	 &&\begin{pmatrix} 
	 1 0 0 1\\ 0 1 1 0
	 \end{pmatrix},\quad
	 &&\begin{pmatrix} 
	 1 0 0 1 \\ 0 1 1 1
	 \end{pmatrix},\quad
	 &&\begin{pmatrix} 
	 1 0 1 0 \\ 0 1 0 0
	 \end{pmatrix},
	 \\[1ex]
	 &\begin{pmatrix} 
	 1 0 1 0 \\ 0 1 0 1
	 \end{pmatrix},
	 &&\begin{pmatrix} 
	 1 0 1 1\\ 0 1 1 0
	 \end{pmatrix},
	 &&\begin{pmatrix} 
	 1 0 1 1\\ 0 1 1 1
     \end{pmatrix},
	 &&\begin{pmatrix} 
	 1 1 0 0\\ 0 0 1 1
	 \end{pmatrix},
	 &&\begin{pmatrix} 
	 1 1 0 1\\ 0 0 1 1
	 \end{pmatrix},
	 \\[1ex]
	 &\begin{pmatrix} 
	 1 0 0 0\\ 0 0 0 1
	 \end{pmatrix},
	 &&\begin{pmatrix}
	 1 0 1 0\\ 0 0 0 1
	 \end{pmatrix},
	 &&\begin{pmatrix} 
	 0 1 0 0\\ 0 0 1 0
	 \end{pmatrix},
	 &&\begin{pmatrix} 
	 0 1 0 1\\ 0 0 1 0
	 \end{pmatrix},
	 &&\begin{pmatrix} 
	 0 0 1 0\\ 0 0 0 1
	 \end{pmatrix}.
      \end{alignat*}

      The lines $\mL_2$ are exactly the 15 following one-dimensional isotropic subspaces:
      \begin{alignat*}{6}
	 & (0001), &\quad& (0100), &\quad& (0111), &\quad& (1010), &\quad& (1101), \\
	 & (0010), &\quad& (0101), &\quad& (1000), &\quad& (1011), &\quad& (1110), \\
	 & (0011), &\quad& (0110), &\quad& (1001), &\quad& (1100), &\quad& (1111).
      \end{alignat*}

      \item 
      Let $n=3$. Then $V=\FF_2^6$ and the sets $\mP_3$ and $\mL_3$ consist of 135 point and $315$ lines, respectively (see the Appendix~\ref{appendix:CR}).
\end{itemize}
\end{example}

\begin{definition}
   For the symplectic dual polar space $\mG_n$ we define its \define{collinearity graph} $\Gamma$ as the graph whose vertices are the points of $\mP_n$. 
   Two vertices are called \define{adjacent} if and only if the corresponding points are collinear in $\mG_n$.
   This is also called the \define{incidence graph} or \define{Menger graph}, see~\cite{cox}.
\end{definition} 
For example, $\mG_1$ consists of one line and three points, so its collinearity graph is a triangle.

\begin{cor}[Number of $(n-1)$-dimensional totally isotropic subspaces]
   \label{cor:NHP}
   The number of different $(n-1)$-dimensional totally isotropic subspaces of a $2n$-dimensional symplectic space $V$ is
   \begin{equation*}
      \#(\text{totally isotropic $(n-1)$-dimensional subspaces of } V) 
      = 
      \frac{(2^n-1)|\mP_n|}{3} 
   \end{equation*}
   where 
   $|\mP_n|$ is the number of all maximal totally isotropic subspaces of $V$ which is calculated in Corollary~\ref{cor:NN}.
\end{cor}

\begin{proof}
   By Corollary~\ref{cor:tres} every $(n-1)$-dimensional isotropic subspace of $V$ is contained in exactly 3 maximal isotropic subspaces.
   Moreover, every maximal isotropic subspace contains exactly $2^n-1$ subspaces of dimension $1$, hence, by passing to the algebraic complement, it has exactly $2^n-1$ subspaces of dimension $n-1$.
   It follows for the collinearity graph $\Gamma$ that every line contains exactly $3$ points and that every point is contained in $2^n-1$ lines.
   Therefore, $\Gamma$ is an incidence structure and, 
   by \cite[Formula (9')]{Dem}, 
   $3 |\mL_n| = (2^n-1)|\mP_n|$ which proves the formula.

\end{proof}

\begin{definition}
   Let $P$ and $Q$ be vertices in a connected graph $\Gamma$.
   A \define{path from $P$ to $Q$ of length $n$} is an ordered set of vertices $V_0=P,\, V_1, \dots, V_n=Q$ such that $V_{i-1}$ and $V_i$ are connected by an edge. 
   The minimal length of all paths connecting $P$ and $Q$ is called the 
   \define{distance} between $P$ and $Q$.
\end{definition}

\begin{definition} 
   \label{def:Gammak}
   Consider the collinearity graph $\Gamma$ associated to the symplectic dual polar space $\mG_n$.
Fix a vertex $X_0$ in the graph $\Gamma$ and consider the induced subgraph $\Gamma_k$ formed by the vertices of $\Gamma$ which have distance $k$ from $X_0$.
\end{definition}

In Theorem \ref{thim} we resume some important properties of the collinearity graph $\Gamma$ associated to the symplectic dual polar space $\mG_n$.
We dedicate the rest of the section to prove them.
Although they seem to be well-known, we could not find easily accessible proofs in the literature. 
For instance, Claim~\ref{item:1} of Theorem \ref{thim} is attributed to the unpublished report~\cite{BardoeIvanov}.
 
\begin{theorem}\label{thim}
   Let $\Gamma$ be the collinearity graph associated to the symplectic dual polar space $\mG_n$ and let $\Gamma_k$ as in Definition~\ref{def:Gammak} for a fixed vertex $X_0$.
   For maximal isotropic subspaces $P$, $Q$ and $R$ of $V$ we have the following properties:
   \begin{enumerate}
      \item\label{item:1}
      $\dim(X_0\cap P) = n-k$ if and only if $P\in\Gamma_k$.

      \item\label{item:2}
      $P$ and $Q$ belong to the same connected component of $\Gamma_k$ if and only if
      $X_0\cap P = X_0\cap Q$.

      \item\label{item:3}
      The induced subgraph $\Gamma_n$ consists of exactly one connected component and 
      $\Gamma_1$ consists of exactly $2^n-1$ disjoint connected components.

      \item\label{item:4}
       Suppose $P,Q, R$ are pairwise different and collinear.
       Then at least two of the spaces belong to the same $\Gamma_k$ and if  $P,Q\in\Gamma_k$, then $R\in\Gamma_{k-1}$.
   \end{enumerate}
\end{theorem}

Claim~\ref{item:1} is shown in Corollary~\ref{cor:components}, 
claim~\ref{item:2} is shown in Theorem~\ref{thm:connectedcomponents},
claim~\ref{item:3} is shown in Corollary~\ref{cor:gamma1}
and
claim~\ref{item:4} is shown in Lemma~\ref{lem:collinear} and Theorem~\ref{thm:collinear}.

\begin{prop}
   \label{cor:pGamma}
   Let $P$ be a maximal totally isotropic subspaces of $V$ with $\dim(X_0\cap P) = n-k$ for some $k\in\{0,1,\dots, n\}$.
   Then $P\in\Gamma_k$, in particular, the collinearity graph $\Gamma$ is connected and the set of its vertices is the disjoint union of the vertices of the $\Gamma_k$.
\end{prop}
\begin{proof}
   According to the proof of Theorem~\ref{thm:SympBasis} we can choose a symplectic basis
   $x_1,\dots, x_n, y_1,\dots, y_n$
   of $V$ and numbers $a_{jm}\in\{0,1\}$ with $a_{jm}=a_{mj}$ such that
   \begin{equation*}
      X_0 = \linspan\{x_1,\dots, x_k\} \oplus \linspan\{ x_{k+1},\dots, x_{n}\},
      \qquad
      P = \linspan\{p_1,\dots, p_k\} \oplus \linspan\{ x_{k+1},\dots, x_{n}\}
   \end{equation*}
   and $p_j = y_j + \sum_{m=1}^k a_{jm}x_m$.
   Let us define the sequence of subspaces $Q_0 := X_0$ and 
   \begin{equation*}
      Q_m := \linspan\{p_1,\dots, p_m\} \oplus \linspan\{ x_{m+1},\dots, x_{n}\},
      \qquad m=1,\dots, k.
   \end{equation*}
   Clearly, all $Q_m$ are maximal totally isotropic subspaces, $\dim(Q_{m-1}\cap Q_{m}) = n-1$ for $m=1,\dots k$ and $X_0\cap P = \bigcap_{m=0}^k Q_m$.

   The subspaces $Q_0,\dots, Q_k$ form a path from $X_0$ to $P$, hence $P\in\Gamma_{k'}$ for some $k'\le k$.
   Now assume that $k'<k$.
   Then we have maximal totally isotropic subspaces $Q_0=X_0,\, Q_1,\, \dots,\, Q_{k'}=p$ with
   $\dim(Q_{m-1}\cap Q_{m})=n-1$ for $m=1,\dots, k'$.
   Now we show by induction that $\dim(X_0\cap Q_1\cap \dots \cap Q_{m})\ge n-m$ for $m=1,\dots, k'$.
   By definition of $Q_1$, we know that 
   $\dim(X_0\cap Q_{1})= n-1$.
   Now assume that the inequality is true for some $1\le m < k'$.
   Then
   \begin{align*}
      \dim(X_0\cap Q_1\cap \dots \cap Q_{m+1})
      & = \dim(X_0\cap Q_1\cap \dots \cap Q_{m}) + \dim Q_{m+1} 
      - \dim( (X_0\cap Q_1\cap \dots \cap Q_m) +  Q_{m+1})
      \\
      & \ge (n-m) + n - \dim(Q_m + Q_{m+1})
      = (n-m) + n - (n+1)
      = n-(m+1).
   \end{align*}
   We obtain
   \begin{equation*}
      n-k = \dim(X_0\cap P)
      \ge \dim(X_0\cap Q_1 \cap \dots \cap Q_{k'-1}\cap P) \ge n-k'
   \end{equation*}
   in contradiction to our assumption $k'<k$.
\end{proof}

\begin{cor}\label{cor:components}
   Let $P$ be a maximal totally isotropic subspace.
   Then $\dim(X_0\cap P) = n-k$ if and only if $P\in\Gamma_k$.
\end{cor}
\begin{proof}
   This follows immediately from Proposition~\ref{cor:pGamma} and the trivial fact that $\Gamma_m\cap \Gamma_\ell = \varnothing$ for $m\neq \ell$.
\end{proof}

\begin{theorem}\label{thm:connectedcomponents}
   Let $P,Q$ be maximal totally isotropic subspaces of $V$.
   Then $P$ and $Q$ belong to the same connected component of $\Gamma_k$ if and only if
   $X_0\cap P = X_0\cap Q$.
\end{theorem}
\begin{proof}
   First let us assume that $X_0\cap P = X_0\cap Q$ and that 
   $\dim(X_0\cap P) = \dim(X_0\cap Q) = n-k$.
   By Theorem~\ref{thm:SympBasis} there exists a symplectic basis $x_1,\dots, x_n,y_1,\dots, y_n$ of $V$ and numbers $a_{jm},\ b_{jm}\in\{0,1\}$ with $a_{jm}=a_{mj}$ and $b_{jm}=b_{mj}$ such that
   \begin{alignat*}{2}
      X_0 &= \linspan\{ x_1,\dots, x_k\} \oplus \linspan\{ x_{k+1},\dots, x_n\} &= \xi \oplus U_0,\\
      P &= \linspan\{ p_1,\dots, p_k\} \oplus \linspan\{ x_{k+1},\dots, x_n\}   &= \pi \oplus U_0,\\
      Q &= \linspan\{ q_1,\dots, q_k\} \oplus \linspan\{ x_{k+1},\dots, x_n\}   &= \eta \oplus U_0,
   \end{alignat*}
   with 
   $\xi= \linspan\{ x_1,\dots, x_k\}$,
   $\pi= \linspan\{ p_1,\dots, p_k\}$,
   $\eta = \linspan\{ q_1,\dots, q_k\}$,
   $U_0 = \linspan\{ x_{k+1},\dots, x_n\}$,
   $p_m = y_m + \sum_{\ell=1}^m a_{m\ell}x_\ell$,
   $q_m = y_m + \sum_{\ell=1}^m b_{m\ell}x_\ell$.

   Let us set 
   $W := \linspan\{ x_{1},\dots, x_k, y_{1},\dots, y_k\}$ and 
   $U := \linspan\{ x_{k+1},\dots, x_n, y_{k+1},\dots, y_n\}$.
   Then clearly $W$ and $U$ are symplectic spaces of dimension $2k$ and $2(n-k)$, respectively, and
   $V=W\operp U$.
   Moreover, $\xi$, $\pi$ and $\eta$ are maximal totally isotropic subspaces of $W$
   and $U_0$ is a maximal totally isotropic subspace of $U$.
   Let $j=\dim(\pi\cap \eta)$.

   As in the proof of Proposition~\ref{cor:pGamma}
   (with $W$, $\eta$ and $\pi$ in place of $V$, $X_0$ and $P$),
   we can find a chain of $k$-dimensional maximal totally isotropic subspaces
   $\rho_0=\eta,\, \rho_1,\dots, \rho_{j-1},\, \rho_j = \pi$ of $W$ such that
   $\dim(\rho_{m-1}\cap \rho_m) = k-1$ for $m=1,\dots, j$
   and $\bigcap_{m=0}^j \rho_m = \eta\cap \pi$.
   Now we define the spaces $R_m := \rho_m\oplus U_0 \subset V$.
   Clearly, all $R_m$ have dimension $n$ and are totally isotropic subspaces of $V$.
   Observe that $R_m\cap X_0 = U_0$, hence $R_m\in\Gamma_k$.
   Moreover, 
   $R_m\cap R_{m'} = (\rho_m\cap \rho_{m'})\oplus U_0$, hence
   $\dim(R_{m-1}\cap R_{m}) = (k-1)+(n-k) = n-1$ for all $m=1,\dots, j$ and 
   $\bigcap_{m=0}^j R_m = Q\cap P$.
   Hence $R_0=Q,\, R_1,\dots, R_j=P$ is a path in $\Gamma_k$ from $Q$ to $P$.
   Actually, Proposition~\ref{cor:pGamma} shows that it is a shortest path between $P$ and $Q$.
   In particular $P$ and $Q$ belong to the same connected component of $\Gamma_k$ and the distance between $P$ and $Q$ within $\Gamma_k$ is $j = \dim(\pi\cap\eta) = \dim(P\cap Q)-\dim(P\cap Q\cap X_0)$.%
   \smallskip

   Now let us assume that $P$ and $Q$ belong to the same connected component of $\Gamma_k$.
   Hence $\dim(X_0\cap P) = \dim(X_0\cap Q)$ by Corollary~\ref{cor:components}.
   \\
   {\em Case 1.}
   Assume that $P$ and $Q$ are collinear, that is, $\dim(P\cap Q)=n-1$.
   Let us suppose that $X_0\cap P \neq X_0\cap Q$.
   Then 
   \begin{align*}
      n-k &= \dim (X_0\cap P) > \dim(X_0\cap P \cap Q)
      = \dim(X_0\cap P) + \dim(Q) - \dim((X_0\cap P) + Q)
      \\
      &\ge (n-k)+n - \dim(P + Q)
      = 2n -k - (n + 1)
      = n-k-1,
   \end{align*}
   hence
   $\dim(X_0\cap P \cap Q) = n-k-1$.
   Now choose a basis $x_{k+2},\dots, x_n$ of $X_0\cap P\cap Q$.
   Next choose $x_k, x_{k+1}\in X_0$ such that
   $x_k\in Q\setminus P$ and 
   $x_{k+1}\in P\setminus Q$.
   Clearly, $x_k,\dots, x_n$ are linearly independent and we can complete them to a symplectic basis 
   $x_1,\dots, x_n, y_1, \dots, y_n$ of $V$ such that 
   \begin{alignat*}{2}
      X_0 &= \linspan\{ x_1,\dots, x_n\},
      \\
      P &= \linspan\{ p_1,\dots, p_k\} \oplus \linspan\{ x_{k+1}\} \oplus \linspan\{ x_{k+2},\dots, x_n\}   
      &= \pi \oplus U_0,\\
      Q &= \linspan\{ q_1,\dots, q_k\} \oplus \linspan\{ x_{k} \} \oplus \linspan\{ x_{k+2},\dots, x_n\} 
      &= \eta \oplus U_0
   \end{alignat*}
   with 
   $\pi = \linspan\{ p_1,\,\dots, p_k, x_{k+1}\}$,
   $\eta = \linspan\{ q_1,\,\dots, q_k, x_{k}\}$,
   $U_0 = \linspan\{ x_{k+2},\dots, x_n\}$,
   and vectors 
   $p_1,\dots, p_k\in\linspan\{x_1,\dots, x_k,y_1,\dots, y_k\}\setminus X_0$ and
   $q_1,\dots, q_k\in\linspan\{x_1,\dots, x_{k-1}, x_{k+1},y_1,\dots, y_{k-1}, y_{k+1}\}\setminus X_0$.

   In particular, we have that
   $p_j, q_j \in\linspan \{x_1,\dots, x_{k+1},y_1,\dots, y_{k+1}\}=: W$.
   By assumption, $\dim (P\cap Q) = n-1$, $\dim U_0 = n-k-1$ and $U_0\subset P\cap Q$.
   Since $P\cap Q = (\pi\cap \eta)\oplus U_0$, it follows that 
   $\dim (\pi\cap \eta) = k+1$.
   Now let us choose a basis $z_1,\dots, z_{k+1}$ of $\pi\cap\eta\subset W$.
   Then there are coefficients $a_{j\ell},\ b_{j\ell}\in\{0,1\}$ such that
   \begin{equation*}
      z_j = \sum_{\ell =1}^{k+1} ( a_{j\ell}x_\ell + b_{j\ell}y_\ell),
      \qquad
      j=1,\dots, k+1.
   \end{equation*}
   Since $x_{k+1}\in P$, $x_k\in Q$ and all $z_j$ belong to the totally isotropic subspaces $P$ and $Q$, we obtain
   \begin{equation*}
      0 = [z_j, x_{k+1}] = b_{j,k+1}
      \quad\text{and}\quad
      0 = [z_j, x_{k}] = b_{jk},
      \qquad j=1,\dots, k+1.
   \end{equation*}
   Therefore
   \begin{equation*}
      z_j = w_j + \sum_{\ell =1}^{k+1} a_{j\ell}x_\ell
      \quad\text{with}\quad
      w_j\in\linspan\{y_1,\dots, y_{k-1}\},
      \qquad
      j=1,\dots, k+1.
   \end{equation*}
   Since $w_1,\dots, w_{k+1}$ are linearly dependent, there are $\mu_1,\dots, \mu_{k+1}$ such that
   $0 = \mu_1 w_1 + \dots + \mu_{k+1}w_{k+1}$ such that not all $\mu_1,\cdots, \mu_{k+1}$ are zero.
   On the other hand, the $z_1,\dots, z_{k+1}$ are linearly independent, so
   \begin{equation*}
      0\neq z_0 := \mu_1 z_1 + \dots + \mu_{k+1}z_{k+1}
      = \sum_{j=1}^{k+1} \sum_{\ell=1}^{k+1} \mu_j a_{j\ell} x_\ell
      \in x_0.
   \end{equation*}
   But then $\dim(\linspan\{z_0\}\oplus U_0) = n-k$ and
   $\linspan\{z_0\}\oplus U_0 \subset X_0\cap Q \cap P$.
   In particular, we obtain
   $\dim(X_0\cap P \cap Q) \ge n-k$, in contradiction to
   $\dim(X_0\cap P \cap Q) = n-k-1$.

   \smallskip

   {\em Case 2.}
   Let $P,\, Q$ be in the same connected component of $\Gamma_k$.
   Then there exists a chain of subspaces 
   $R_0=Q,\, R_1, \dots, R_j=P$ in the connected component of $\Gamma_k$ with 
   $\dim(R_{m-1}\cap R_m) = n-1$ for all $m=1,\dots, j$.
   By Case~1, it follows that 
   $Q\cap X_0 = R_1\cap X_0 = \dots  = R_{j-1}\cap X_0 = P\cap X_0$.
\end{proof}

\begin{cor}\label{cor:gamma1}
   $\Gamma_n$ consists of exactly one connected component.
   $\Gamma_1$ consists of exactly $2^n-1$ components.
\end{cor}
\begin{proof}
   Recall that $P\in \Gamma_n$ if and only if $\dim(X_0\cap P)=0$.
   Hence $P\cap X_0 = Q\cap X_0 = \{0\}$ for all $P,Q\in\Gamma_n$, hence they belong to the same connected component.

   The number of connected components of $\Gamma_1$ is equal to the different $(n-1)$-dimensional subspaces of $X_0$, which in turn is equal to the number of different one-dimensional subspaces of $X_0$ (by taking algebraic complements), which is equal to the number of non-zero elements of $X_0$ which is equal to $2^n-1$.
\end{proof}

Next we will prove the so-called $n$-gon property. 
First we observe the following.
\begin{lemma}
   \label{lem:collinear}
   Let $P,Q, R$ collinear.
   Then at least two of the spaces belong to the same $\Gamma_k$.
\end{lemma}
\begin{proof}
   If two of the spaces are equal, then the claim is clear.
   Now assume that they are three different subspaces and that
   $P\in\Gamma_k$, $Q\in\Gamma_\ell$ and $R\in\Gamma_m$.
   Without restriction we assume that $m=\min\{k,\ell, m\}$.
   Since $P,Q, R$ are collinear, they are at distance $1$ from each another.
   Since $k$ is the distance of $P$ from $X_0$, it follows that
   $m\le k \le m+1$. Analogously, $m\le \ell \le m+1$.
\end{proof}

\begin{theorem}\label{thm:collinear}
   Let $P,Q, R$ be collinear and pairwise different.
   If $P,Q\in\Gamma_k$, then $R\in\Gamma_{k-1}$.
\end{theorem}
\begin{proof}
   Since $P,Q, R$ are collinear, it follows that
   $P\cap Q = Q\cap R = Q \cap R = P\cap Q\cap R =: U$
   and $\dim(U) = n-1$.
   Moreover, $P,Q$ belong to the same connected component of $\Gamma_k$, hence
   $\dim(X_0\cap P) = \dim(X_0\cap Q) = n-k$
   and
   $X_0\cap P = X_0\cap Q$ and therefore
   \begin{equation*}
      X_0\cap P = X_0 \cap Q = X_0\cap P \cap Q
      = X_0\cap U
      \subseteq
      X_0\cap R.
   \end{equation*}
   Let $x_1,\dots, x_n, y_1,\dots, y_n$ be a symplectic basis of $V$ such that
   $x_{k+1},\dots, x_n$ is basis of $X_0\cap P\cap Q$ and $x_1,\dots, x_n$ is basis of $X_0$.
   Let $W=\linspan\{x_1,\dots, x_k, y_1,\dots, y_k\}$.
   Then, by Theorem~\ref{thm:SympBasis}, there are $k$-dimensional spaces $\widetilde P, \widetilde Q\subseteq W\setminus X_0$ with $\dim(\widetilde P \cap \widetilde Q)=k-1$ such that
   \begin{equation*}
      P = \widetilde P\oplus (X_0\cap P),\qquad
      Q = \widetilde Q\oplus (X_0\cap Q).
   \end{equation*}
   Let us choose a basis $z_1,\dots, z_{k-1}$ of $\widetilde P\cap \widetilde Q$; without restriction, we may assume that it is of the form
   \begin{equation}
      \label{eq:z}
      z_j = \widetilde y_j + \sum_{\ell=1}^{k} a_{j\ell} x_\ell,
      \qquad j=1,\dots, k-1,
   \end{equation}
   with $\widetilde y_j,\, j=1,\dots, k$ being a basis of $\linspan \{ y_1,\dots, y_k\}$.
   Observe that this means that
   \begin{equation*}
      U = P\cap Q = \linspan\{z_1,\dots, z_{k-1}, x_{k+1},\dots, x_n\}.
   \end{equation*}
   Now we choose $z_k^P\in W$ to complete $\{z_1,\dots, z_{k-1}\}$ to a basis of $\widetilde P$.
   Note that $z_k^P\notin X_0$, because otherwise $\dim(\widetilde P \cap X_0)\ge 1$, hence $\dim(P\cap X_0)\ge (n-k)+1$.
   Thus, by \eqref{eq:z} we may assume that
   $z_k^P = \widetilde y_k + \sum_{\ell=1}^k \beta_\ell x_\ell$.
   Analogously, we can choose
   $z_k^Q = \widetilde y_k + \sum_{\ell=1}^k \gamma_\ell x_\ell$
   such that $\{z_1,\dots, z_{k-1}, z_k^Q\}$ is a basis of $\widetilde Q$.
   In summary, we have now that
   \begin{equation*}
      P = U \oplus\linspan\{z_k^P\},\qquad
      Q = U \oplus\linspan\{z_k^Q\}
   \end{equation*}
   are two different maximal totally isotropic subspaces containing the $(n-1)$-dimensional isotropic subspace $U$.
   Since $R$ also contains $U$, by \eqref{eq:3Lagrangians} it must be of the form
   \begin{equation*}
      R = U \oplus\linspan\{z_k^P+z_k^Q\}
      = \linspan\{z_1, \dots, z_{k-1}\}\oplus\linspan\{z_k^P+z_k^Q\}\oplus\linspan\{x_{k+1},\dots, x_n\}.
   \end{equation*}
   The $z_1,\dots, z_{k-1}$ are linearly independent over $X_0$ and 
   $z_k^P + z_k^Q = \sum_{\ell=1}^k (\beta_\ell + \gamma_\ell) x_\ell \in X_0$,
   so it follows that $\dim(R\cap X_0)=n-k+1$, hence $R\in\Gamma_{k-1}$.
   (It is sufficient to know that $z_k^P + z_k^Q\in X_0$ because this implies that $\dim(R\cap X_0)\ge n-k+1$. On the other hand $\dim(R\cap X_0)\le n-k+1$ because $R$ has distance $1$ from $P$.)
\end{proof}

\begin{remark}
   For every $k=0, \dots, n-1$ and every $R\in\Gamma_{k}$, there exist $P,Q$ in $\Gamma_{k+1}$ such that $R,P,Q$ are collinear.
\end{remark}
\begin{proof}
   By Corollary~\ref{cor:tres} and Theorem~\ref{thm:collinear} it is suffices to find $P\in\Gamma_k$ which is collinear with $R$.
   As in the proof of Theorem~\ref{thm:connectedcomponents}, we can choose a symplectic basis $x_1,\dots, x_n,y_1,\dots, y_n$ such that 
   \begin{alignat*}{2}
      X_0 &= \linspan\{ x_1,\dots, x_k\} \oplus \linspan\{ x_{k+1},\dots, x_n\},
      \\
      R &= \linspan\{ r_1,\dots, r_k\} \oplus \linspan\{ x_{k+1},\dots, x_n\}
   \end{alignat*}
   with 
   $r_m = y_m + \sum_{\ell=1}^m a_{m\ell}x_\ell$ with $a_{m\ell}\in\FF_2$.
   Then clearly
   \begin{equation*}
      P = \linspan\{ r_1,\dots, r_{k}\} \oplus \linspan\{ y_{k+1}\} \oplus \linspan\{ x_{k+2},\dots, x_n\}
   \end{equation*}
   has the desired properties.
\end{proof}
   

\section{The universal embedding of $\mG_n$} 
\label{sec:embedding}
\begin{definition}
   An \define{embedding of $\mG_n$} is an $\FF_2$-vector space $E$ together with a map $\theta:\mP_n\to E$ such that
   \begin{enumerate}
      \addtolength{\itemsep}{-.4\baselineskip}
      \item $\theta(P)\neq 0$ for every $P\in\mP_n$,
      \item $E = \linspan\{ \range(\theta) \}$ where $\range(\theta)=\theta(\mP_n)$ is the range of $\theta$.
      \item $\theta(P) + \theta(Q) +\theta(R) = 0$ for every line $L=\{P,Q,R\}\in\mL_n$.
   \end{enumerate}
\end{definition}

Such an embedding can be constructed as follows.
Let $\FF_2(\mL_n)$, $\FF_2(\mP_n)$ be the $\FF_2$-vector spaces freely generated by the lines $\mL_n$ and the points $\mP_n$, respectively.
Since every line $L\in\mL_n$ can be written as $L=\{P,Q, R\}$ where $P,Q, R\in\mP_n$ are the three points contained in $L$, we have the following map
\begin{equation*}
   \eta:\FF_2(\mL_n)\To \FF_2(\mP_n)\,,
   \qquad L=\{P,Q,R\} \mapsto P+Q+R.
\end{equation*}
The quotient $\UGn := \FF_2(\mP_n)/\eta(\FF_2(\mL_n))$ is called the 
\define{universal embedding module} and we define the canonical map
\begin{equation*}
   \theta: \mP_n\To \UGn \,.
\end{equation*}
Clearly this is an embedding of $\mG_n$; it is called its \define{universal embedding}.
Note that any other embedding is a quotient of the universal embedding.
The \define{dimension of the universal embedding of the polar dual space} is $\dim(\UGn)$.

Brouwer proved in 1990 that 
\begin{equation}\label{eqr1}  
   \dim(\UGn)\geq \frac{(2^{n}+1)(2^{n-1}+1)}{3}
\end{equation}
and conjectured that \eqref{eqr1} is in reality an equality.
The conjecture was proved by Li \cite{li} and independently by Blokhuis and Brouwer \cite{brouwer}.

\begin{theorem}[\cite{li, brouwer}]
   \label{BL}
   The dimension of the universal embedding of the polar dual space, $\dim U(\mG_n)$, is
   \begin{equation*}
      \dim U(\mG_n)=\frac{(2^{n}+1)(2^{n-1}+1)}{3}\,.
   \end{equation*}
\end{theorem}

In the rest of this section we describe briefly how Li in \cite{li} used certain vector spaces to count $\dim U(\mG_n)$ and thereby proved Theorem~\ref{BL}.
Since \eqref{eqr1} was already known, 
only the reverse inequality $\le$ had to be proved.

To this end, Li considers the collinearity graph $\Gamma$ defined by $\mG_n$.
As before, we fix a point $X_0$ in $\Gamma$ and we set $\Gamma_k$ to be the set of all points in $\Gamma$ which have distance $k$ from $X_0$.
Since every triangle in $\Gamma$ contains two elements from $\Gamma_k$ and one from $\Gamma_{k-1}$ for some $k=1,\dots, n$, if follows that $\theta(Y) \in\linspan\{\theta(\Gamma_{k})\}$ for every $Y\in\Gamma_{k-1}$.
Thus we have the following filtration of $\UGn$
\begin{equation*}
   \{0\} 
   \subset \linspan\{ \theta(\Gamma_0) \} = \linspan\{ \theta(X_0) \}
   \subset \linspan\{ \theta(\Gamma_1) \}
   \subset \dots
   \subset \linspan\{ \theta(\Gamma_n) \} = \UGn
\end{equation*}
and consequently
\begin{equation*}
   U(\mG_n) \cong \linspan\{ \theta(\Gamma_0) \}
   \oplus \big( \linspan\{ \theta(\Gamma_1) \}/ \linspan\{ \theta(\Gamma_0) \} \big)
   \oplus \cdots
   \oplus \big( \linspan\{ \theta(\Gamma_n) \}/ \linspan\{ \theta(\Gamma_{n-1}) \} \big).
\end{equation*}

Recall that two points $P,Q$ belong to the same connected component of $\Gamma_k$ if and only if $P\cap X_0 = Q\cap X_0$.
Clearly, this is the case if and only if $\theta(P) \equiv \theta(Q)\mod \linspan\{\theta(\Gamma_{k-1})\}$.

Now let $n\ge 3$ and $2\le k \le n-1$ and let $L,M\subset X_0$ be subspaces with $\dim L = n-k-1$ and $\dim M = n-k+2$.
Then there are exactly $7$ subspaces $L\subset R_j\subset M$ with $\dim R_j=n-k$, and
$\sum_{j=1}^7 \theta (\widetilde R_j) \equiv 0 \mod\linspan\{\theta(\Gamma_{k-1})\}$
where $\widetilde R_j$ is any maximal totally isotropic subspace of $V$ with $\widetilde R_j\cap X_0 = R_j$. 

For $1\le i \le n$ we set $\mathcal W_i$ to be the $\FF_2$-vector space freely generated by all $i$-dimensional subspaces of $X_0$ and 
for $1\le i < j \le n$ we set $\mathcal W_{ij}$ to be the $\FF_2$-vector space freely generated by all flags $X<Y$ in $X_0$ where $\dim X = i$ and $\dim Y = j$.
Let $\{e_L\}$ be the natural basis of $\mathcal W_i$ and $\{e_{X<Y}\}$ the natural basis of $\mathcal W_{ij}$.
Let us define the incidence map
\begin{equation*}
   \phi_{n-k}: \mathcal W_{n-k-1, n-k+2}\to \mathcal W_{n-k},\qquad
   \phi_{n-k}(e_{X<Y}) = \sum_{\substack{X\subset L \subset Y\\ \dim L = n-k}}e_{L}.
\end{equation*}
Moreover, we have a natural surjection
\begin{equation*}
    h_{n-k}: \mathcal W_{n-k} \to \linspan\{\theta(\Gamma_k)\}/\linspan\{\theta(\Gamma_{k-1})\}.
\end{equation*}
It follows from the above that 
$ h_{n-k}\circ \phi_{n-k} = 0$, hence the induced map
\begin{equation*}
   \widetilde h_{n-k}: \mathcal W_{n-k}/\range(\phi_{n-k}) \to \linspan\{\theta(\Gamma_k)\}/\linspan\{\theta(\Gamma_{k-1})\}
\end{equation*}
is well-defined and surjective.
Therefore
\begin{align}
   \nonumber
   \dim \UGn &= 
   \dim\Big( \linspan\{\theta(\Gamma_0)\}\}\Big)
   + \dim\Big( \linspan\{\theta(\Gamma_1)\}/\linspan\{\theta(\Gamma_{0})\}\Big)
   \\
   \nonumber
   &\phantom{= } +
   \sum_{k=2}^{n-1} \dim\Big( 
   \linspan\{\theta(\Gamma_k)\}/\linspan\{\theta(\Gamma_{k-1})\}
   \Big)
   + \dim\Big( \linspan\{\theta(\Gamma_n)\}/\linspan\{\theta(\Gamma_{n-1})\}\Big)
   \\
   \nonumber
   &\le 
   1 + (2^n-1) + \sum_{j=1}^{n-2} \dim\Big( \mathcal W_j/\range(\phi_j) \Big) + 1
   =
   1 + 2^n + \sum_{j=1}^{n-2} \dim( \mathcal W_j) - \dim(\range(\phi_j))
   \\
   \label{eq:E1}
   &=
   \sum_{j=0}^{n} \dim( \mathcal W_j)
   - \sum_{j=1}^{n-2} \dim(\range(\phi_j)),
\end{align}
where in the last step we have used that 
$\dim \mathcal W_0 = \dim \mathcal W_n = 1$,
$\dim \mathcal W_{n-1} = 2^n - 1$.

In order to evaluate the right hand side, Li introduces the following order on $X_0\cong \FF_2^n$.
Let $e_1,\dots, e_n$ be the standard basis of $\FF_2^n$ and let
$v = \alpha_1 e_1 + \dots + \alpha_n e_n$ and 
$w = \beta_1 e_1 + \dots + \beta_n e_n$ in $\FF_2^n$.
Then we define the support and weight of $v$ as
\begin{equation*}
   \supp v = \{j : \alpha_j\neq 0\},
   \qquad
   \wt v = |\supp v|
\end{equation*}
and we set
$m(v) = \min(\supp v)$ and $M(v) = \max(\supp v)$.
We obtain a total order on $\FF_2^n$ by setting $v\succ w$ if and only if there is a $j=1,\dots, n$ such that $\alpha_k=\beta_k$ for all $1\le k < j$ and $(\alpha_j,\, \beta_j) = (1,\, 0)$.

If $L\subseteq \FF_2^n$ is a subspace, we set 
$\supp L  = \bigcup_{v\in L} \supp v$ and
$m(L)  = \{ m(v) : v\in L\}$.
It is not hard to see that $\dim L = | m(L) |$.
The \define{reduced echelon basis of $L$} is the unique basis $v_1,\dots, v_k$ such that $m(v_j)$ is strictly increasing and $m(v_j)\neq \supp\{v_1,\,\dots,\, v_{j-1},\, v_{j+1},\, \dots,\, v_k\}$.
This basis is obtained easily if we take an arbitrary basis of $L$, form the matrix whose rows are these basis vectors and apply the Gau\ss-Jordan procedure to obtain a reduced row-echelon matrix. The rows of this new matrix form the reduced echelon basis of $L$.
Now we define an order on the subspaces of $\FF_2^n$ as follows. 
Let $L, L'$ be subspaces of $\FF_2^n$ with $\dim L = \dim L' = k$ and reduced echelon basis
$v_1,\dots, v_n$ and $v'_1,\dots, v'_n$ respectively.
Then we say that $L \succ L'$ if there is $j=1,\dots, n$ such that $v_k=v'_k$ for $k>j$ and $v_j\succ v'_j$.

Recall that an element $\Delta\in\range\phi_k$ is a formal sum of $k$-dimensional subspaces of $X_0$.
Let us set $A_k = \{\max \Delta : \Delta \in\range(\phi_k)\}$.
Then it is not hard to see that $\dim(\range(\phi_k)) = |A_k|$.

Now we define $\ma{N}^n := \{ L\subset X_0\} \setminus \bigcup_{k=1}^{n-2} A_k$ to be the set of subspaces of $X_0$ which belong to no $A_k$.
Clearly, all the sets $A_k$ are disjoint.
Recall that $\dim \mathcal W_j$ is the number of all subspaces of $X_0$ of dimension $j$.
So we obtain from \eqref{eq:E1}
\begin{equation}
   \label{eq:clever}
   \dim \UGn \le 
   \sum_{j=0}^{n} \dim( \mathcal W_j) - \sum_{j=0}^{n} |A_j|
   =  | \ma{N}^n|.
\end{equation}

Li gives a clever description of the elements in $\ma{N}^n$ using the reduced echelon basis as follows,
see also~\cite{mcclurg}.

\begin{theorem}
   Let $L$ be a $k$-dimensional subset of $X_0$ with reduced echelon basis $v_1 \succ  \dots \succ  v_k$.
   Then $L\in \ma{N}^n$ if and only if the following four conditions are satisfied:
   \begin{itemize}
      \item[(N1)]
      $\wt v_j \le 2$ for all $j=1,\dots, k$.

      \item[(N2)]
      If $v_r \succ  v_s$ and 
      $\wt v_r  = \wt v_s = 2$, then $M(v_r) \le M(v_s)$.

      \item[(N3)]
      If $v_r \succ  v_s \succ  v_t$,
      $\wt v_r  = \wt v_s  = \wt v_t = 2$ and $M(v_r) = M(v_s) < M(v_t)$,
      then $m(v_t) > M(v_r)$.

      \item[(N4)]
      If $v_r \succ  v_s \succ  v_t \succ  v_u$ and 
      $\wt v_r  = \wt v_s  = \wt v_t = \wt v_u = 2$,
      then it is impossible that $M(v_r) = M(v_s) = M(v_t) < m(v_u)$. 
   \end{itemize}
\end{theorem}
\noindent
Note that the last condition in (N4) is equivalent to $M(v_r) = M(v_s) = M(v_t) < M(v_u)$ by condition (N3).
Then Li shows how $\ma{N}^{n+1}$ can be constructed from $\ma{N}^n$ and thus is able to show that  
\begin{equation*}
   |\ma{N}^n| = \frac{(2^{n}+1)(2^{n-1}+1)}{3}\,,
\end{equation*}
proving the formula in Theorem~\ref{BL}. 
In addition, it follows that the functions $\widetilde h_{n-k}$ are bijections and that
$\range(\phi_k) = \ker( h_{n-k})$.
\smallskip

If we modify slightly Li's construction of $\ma{N}^{n+1}$ from $\ma{N}^n$, then it is analogous to how we constructed $W^{n+1}$ from $W^n$.
In the next section we will show how this allows us to construct a bijection between $W^{n+1}$ and $\ma{N}^n$.

\section{Bijection between words and vector spaces} 
\label{sec:bijection}
As in \cite{li} we set $g(n) := |\ma{N}^n|$.
Before we continue, let us give some examples of $\ma{N}^n$.
We use the notation of Example~\ref{ex:VS}.
\begin{itemize}
   \item The set of all $L\in \ma{N}^1$ is $(0)$ and $(1)$ and $g(1) = |\ma{N}^1| = 2 = g_W(2)$.

   \item The set of all $L\in \ma{N}^2$ is, in ascending order, $(00)$, $(01)$, $(10)$, $(11)$, 
   $\left(\begin{smallmatrix}10\\ 01 \end{smallmatrix}\right)$
   and $g(2) = |\ma{N}^2| = 5 = g_W(3)$.

   \item The set of all $L\in \ma{N}^3$ is, in ascending order, 
   \begin{align*}
   &(000),\quad (001),\quad (010),\quad (100),\quad (011),\quad (101),\quad (110), \\
   &\begin{pmatrix}010\\ 001 \end{pmatrix},\quad
   \begin{pmatrix}100\\ 001 \end{pmatrix},\quad
   \begin{pmatrix}110\\ 001 \end{pmatrix},\quad
   \begin{pmatrix}100\\ 010 \end{pmatrix},\quad
   \begin{pmatrix}101\\ 010 \end{pmatrix},\quad
   \begin{pmatrix}100\\ 011 \end{pmatrix},\quad
   \begin{pmatrix}101\\ 011 \end{pmatrix},\quad
   \begin{pmatrix}100\\ 010 \\ 001 \end{pmatrix}
   \end{align*}
   and $g(3) = |\ma{N}^3| = 15 = g_W(4)$.
\end{itemize}

Now we show how $\ma{N}^{n}$ can be constructed from $\ma{N}^{n-1}$ for $n\ge 2$.
For this to be analogous to the process for the passage from $W^n$ to $W^{n+1}$, we have to modify Li's procedure slightly.

In what follows, we identify vector spaces $L$ with the matrix $A$ whose rows consist of the reduced echelon basis of $L$ and we use the following notation:
If $v\in \FF_2^{n-1}$, then we denote by 
$\widetilde v\in\FF_2^n$ the vector obtained from $v$ by appending a $0$ 
and by 
$\widehat v\in\FF_2^n$ the vector obtained from $v$ by inserting a $0$ between the last and second to last component of $v$.
The $k$th unit vector in $\FF_2^n$ is denoted by $e_k^n$.

We will say that a vector space $\widetilde L\in \ma{N}^{n}$ \define{is in Case~$k'$} for $k'=1,\,\dots,\, 8$, if it is constructed from a vector space $L\in \ma{N}^{n-1}$ as described in the Cases~$k'$ below.
For examples of these constructions, see the ones in Example~\ref{ex:Psi}.

\begin{itemize}
   \item{\bf Case $1'$.}
   Take an arbitrary vector space in $\ma{N}^{n-1}$ with reduced echelon basis $v_1\succ \dots \succ v_{k}$. 
   Append $0$ to each of these vectors in order to obtain $\widetilde v_1\succ \dots \succ \widetilde v_{k}$.
   Then clearly $\widetilde L\in \ma{N}^n$ and $\dim L = \dim \widetilde L$.
   We denote this construction by
   \begin{equation*}
      \alpha^{n-1}_{n,0} :\ \ma{N}^{n-1}\to \ma{N}^{n},
      \qquad L\mapsto \widetilde L.
   \end{equation*}
   Note that each vector space $\widetilde L$ obtained in this way has the form 
   $\displaystyle
   \begin{pmatrix}
      & &  0 \\
      & L  & \vdots \\
      & &  0 
   \end{pmatrix}
   $
   for some vector space $L\in \ma{N}^{n-1}$.
   The total number of such vector spaces $\widetilde L$ is $g(n-1)$.

   \item{\bf Case $2'$.}
   Take an arbitrary vector space in $\ma{N}^{n-1}$ with reduced echelon basis $v_1\succ \dots \succ v_{k}$. 
   Append $0$ to each of these vectors in order to obtain $\widetilde v_1\succ \dots \succ \widetilde v_{k}$ and augment this basis by $e_{n}^n$ to a reduced echelon basis of $\widetilde L:=\linspan\{\widetilde v_1,\dots, \widetilde v_{k}, e_{n}^n\}$.
   Then clearly $\widetilde L\in \ma{N}^n$ and $\dim \widetilde L = \dim L + 1$.
   We denote this construction by
   \begin{equation*}
      \alpha^{n-1}_{n,1} :\ \ma{N}^{n-1}\to \ma{N}^{n},
      \qquad L\mapsto \widetilde L.
   \end{equation*}
   Note that each vector space $\widetilde L$ obtained like this has the form 
   $\displaystyle
   \begin{pmatrix}
      & & & 0 \\
      & L & & \vdots \\
      &   & & 0 \\
      0 & \dots & 0 & 1
   \end{pmatrix}
   $
   for some vector space $L\in \ma{N}^{n-1}$.
   The total number of such vector spaces $\widetilde L$ is $g(n-1)$.
   \end{itemize}

   \noindent
   By construction, every vector space obtained in Cases~1' and 2' has either only zeros in the last column or its last line is the vector $e_n^n$.
   \smallskip

   \noindent
   For the remaining cases 3', 4', 5', 6' and 7' we take an arbitrary vector space $L\in \ma{N}^{n-1}$ with reduced echelon basis $v_1\succ \dots \succ v_{k}$ such that $n-1\in\supp L$ and $v_{k}\neq e_{n-1}^{n-1}$.
   This means that the matrix consisting of the row vectors $v_1,\dots, v_{k}$ has at least one $1$ in its last column and the last row is not equal to $e_{n-1}^{n-1}$ and $L$ cannot have been obtained from $(0)$ or $(1)\in \ma{N}^{1}$ by using only the constructions described in Cases~$1'$ and $2'$.

   \begin{itemize}
   \item{\bf Case $3'$.}
   Insert a $0$ in front of the last coordinate of each of the basis vectors in order to obtain $\widehat v_1\succ \dots \succ \widehat v_{k}$ and set 
   $\widetilde L:=\linspan\{\widehat v_1,\,\dots,\, \widehat v_{k}\}$.
   We denote this construction by
   \begin{equation*}
      \alpha^{n-1}_{n-1,0}: L\mapsto \widetilde L.
   \end{equation*}
   Clearly $\widetilde L\in \ma{N}^n$ and $\dim \widetilde L=\dim L$.

   \item{\bf Case $4'$.}
   As in Case $3'$, insert a $0$ in front of the last coordinate of each of the basis vectors in order to obtain 
   $\widehat v_1\succ \dots \succ \widehat v_{k}$ and set $\widetilde L:=\linspan\{\widehat v_1,\,\dots,\, \widehat v_{k},\, e_{n-1}^n\}$.
   We denote this construction by
   \begin{equation*}
      \alpha^{n-1}_{n-1,1}: L\mapsto \widetilde L.
   \end{equation*}
   Clearly $\widetilde L\in \ma{N}^n$ and $\dim \widetilde L=\dim L+1$.

   \item{\bf Case $5'$.}
   As in Case $3'$, insert a $0$ in front of the last coordinate of each of the basis vectors in order to obtain $\widehat v_1\succ \dots \succ \widehat v_{k}$ and set $\widetilde L:=\linspan\{\widehat v_1,\,\dots,\, \widehat v_{k},\, e_{n-1}^n+e_n^n\}$.
   We denote this construction by
   \begin{equation*}
      \alpha^{n-1}_{n-1,2}: L\mapsto \widetilde L.
   \end{equation*}
   Clearly $\widetilde L\in \ma{N}^n$ and $\dim \widetilde L=\dim L+1$.

\end{itemize}

\noindent
The number of vector spaces in each of the Cases $3', 4', 5'$ is 
\begin{align*}
   g(n-1)-\#(\text{$L\in \ma{N}^{n-1}$ with either last column zero or $v_{k}=e_{n-1}^{n-1}$ } )
   = g(n-1) - 2 g(n-2),
\end{align*}
because
$ \#(\text{$L\in \ma{N}^{n-1}$ with last column zero} )
= \#(\text{$L\in \ma{N}^{n-1}$ with $v_{k}=e_{n-1}^{n-1}$} )
= g(n-2)$ as in Case~$1'$.

\begin{itemize}
   \item{\bf Case $6'$.}   Let $L\in \ma{N}^{n-1}$ with reduced echelon basis $v_1\succ \dots \succ v_{k}$ such that either the last column of $A$ has at least two ones (this corresponds to case~6 in \cite{li}), or such that the last column has exactly one 1 and that this 1 is not in the first row (this is a subset of the cases~7b, c and d of \cite{li}).

   \begin{itemize}
      \item{\bf Case $6'$a.}
      Assume that the last column of $A$ has at least two ones. 
      Then every row with a $1$ in its last column must have weight $2$.
      Set $j=\min\{\ell : M(v_\ell) = n-1\}=$ highest row of $A$ with a $1$ in the last column.
      Let $b<n-1$ such that $v_j=e_b^{n-1}+e_{n-1}^{n-1}$.
      For $\ell\neq j$ we let $\widehat v_\ell$ be the vector in $\FF^n$ which is obtained from $v_\ell$ by inserting a $0$ between the last and the second to last component of $v_\ell$.
      Then we set 
      $\widetilde L = \linspan\{ \widehat v_1,\,\dots,\, \widehat v_{j-1},\, e_b^{n-1}+e_{n-1}^{n-1},\, \widehat v_{j+1},\,\dots,\, \widehat v_k\}$.
      In words: we append to $A$ a zero column and then we push the $1$s in the $(n-1)$th column below the $j$th row out to the new $n$th column.
      The matrix $A'$ corresponding to $\widetilde L$ has exactly one $1$ in the second to last column; this occurs in a row with weight $2$, different from the last row and there is at least one row of the form $e_a^n+e_{n}^n$ with $a>j$.
      Note that $\dim \widetilde L = \dim L$. This is the case~6 in \cite{li}.

      \item{\bf Case $6'$b.}
      Assume that the last column of $A$ has at exactly one $1$ and let $j>2$ such that $v_j = e_a^{n-1} + e_{n-1}^{n-1}$.
      Let $b = \max\{ M(v_\ell) : \ell \neq j \}$.

      If $b<a$, then necessarily $j=k$ and we set
      $\widetilde L=\linspan\{\widehat v_1,\,\dots,\, \widehat v_{k-1},\, e_a^{n},\, e_{n-1}^n+e_n^n\}$
      (this is part of case~7a in \cite{li}).
      Note that $\dim \widehat L = \dim L + 1$.

      If $a<b<n-1$ and there exists $\ell > j$ with $b\in\supp v_\ell$, then $v_\ell = v_k =e_b^{n-1}$ and $b\not\in \supp v_m$ for $m\neq\ell$.
      We set
      $\widetilde L = \linspan\{ \widehat v_1,\,\dots,\, \widehat v_{j-1},\, e_a^n+e_b^n,\, \widehat v_{j+1},\,\dots,\, \widehat v_{k-1}, e_{n-1}^n+e_n^n\}$
      (this is part of case~7b in \cite{li}).
      Note that $\dim \widetilde L = \dim L$.

      If $a<b<n-1$ and there exists $\ell < j$ with $b\in\supp v_\ell$, then $\wt v_\ell = 2$ and there is $c<a<b$ such that $v_\ell = e_c^{n-1}+e_b^{n-1}$.
      We set
      $\widetilde L = \linspan\{ \widehat v_1,\,\dots,\, \widehat v_{j-1},\, e_a^n+e_b^n,\, \widehat v_{j+1},\,\dots,\, \widehat v_{k-1}, \widehat v_k,\, e_{n-1}^n+e_n^n\}$
      (this is case~7c in \cite{li}).
      Note that $\dim \widetilde L = \dim L+1$.
   \end{itemize}
   We denote these constructions by 
   \begin{equation*}
      \alpha^{n-1}_{n-1,3} : L\mapsto \widetilde L.
   \end{equation*}

   \item{\bf Case $7'$.}
   Let $L\in \ma{N}^{n-1}$ with reduced echelon basis $v_1\succ \dots \succ v_{k}$ such that $L$ is not in case~6'.
   Then $\wt v_1 = 2$, $n-1\in\supp v_1$ and $n-1\not\in\supp v_j$ for $j\ge 2$.
   That implies that $\wt v_j = 1$ for $j\ge 2$.
   Moreover, it implies that $L$ originates from a vector space $\widehat L\in \ma{N}^{\widehat n}$ in Case~8 using only the constructions described in Cases~$3'$ and $4'$. 
   
   Let $v_1 = e_b^{n-1} + e_{n-1}^{n-1}$.
   If $\dim L = 1$, we set $\widetilde L=\linspan\{e_b^n,\, e_{n-1}^n+e_n^n\}$
   (this corresponds to 1-dimensional vector spaces of case~7a in \cite{li}).

   If $\dim L > 1$ and $v_k=e_a^{n-1}$, we set $\widetilde L=\linspan\{e_b^n + e_a^n, \widehat v_2,\, \dots,\, \widehat v_{k-1},\,  e_{n-1}^n+e_n^n\}$
   (this corresponds to some of the vector spaces of case~7b in \cite{li}).

   We denote these constructions by 
   \begin{equation*}
      \alpha^{n-1}_{n,3} : L\mapsto \widetilde L.
   \end{equation*}

\end{itemize}

\noindent
The total number of vector spaces in the Cases~$6'$ and $7'$ together is
\begin{align*}
   g(n-1)-\#(\text{$L\in \ma{N}^{n-1}$ with either last column zero or $v_{k}=e_{n-1}^{n-1}$ } )
   = g(n-1) - 2 g(n-2).
\end{align*}

\begin{itemize}
   \item{\bf Case $8'$.}
   Let $\widetilde L=(0\cdots 011)\in \ma{N}^n$.
   Clearly this vector space is not contained in the spaces constructed so far.
\end{itemize}

It is not hard to see that the vector spaces constructed above are all pairwise disjoint and that they all belong to $\ma{N}^n$.
Moreover, it can be seen that we obtain every $\widetilde L\in \ma{N}^n$ in exactly one way.

It is clear that the cases $k$ for words and $k'$ for vector spaces correspond to each other.
E.g., appending a $0$ to a given word corresponds to appending a zero column to a vector space (Case 1 and $1'$);
appending a $1$ to a given word corresponds to appending a zero column to a vector space and adding the base vector $x_n$ (Case 2 and $2'$); etc.

Therefore we obtain the following theorem in analogy to Theorem~\ref{thm:wordsconstuction}:

\begin{theorem}
   \label{thm:VSconstuction}
   Let $n\ge 2$.
   Then for every vector space $L\in \ma{N}^{n}$ exactly one of the following holds.
   \begin{enumerate}
      \item 
      There is exactly one sequence of maps $B^2, \dots, B^{n-1}$ such that 
      $L=B^{n-1} \cdots B^{2}(0)$ or $L=B^{n-1} \cdots B^{2}(1)$
      where the $B^j$ are maps of type $\alpha^j_{\ell, a}$ as in the cases above.

      \item 
      There is exactly one $k\le n$ and exactly one sequence of maps $B^k, B^{k+1}, \dots, B^{n-1}$ such that 
      $L=B^{n-1} \cdots B^{k}(\widehat L)$
      where the $B^j$ are maps of type $\alpha^j_{\ell, a}$ as in the cases above and $\widehat L=(0\dots 011)\in \ma{N}^k$.
   \end{enumerate}
\end{theorem}

Now Theorem~\ref{thm:wordsconstuction} and 
Theorem~\ref{thm:VSconstuction} together give a bijection between $W^{n+1}$ and $\ma{N}^n$.

\begin{theorem}
   \label{thm:bijection}
   Let $n\in\NN$.
   We have the following bijection
   \begin{equation*}
      \Psi: W^{n+1}\to \ma{N}^n
   \end{equation*}
   defined as follows.
   \begin{enumerate}
      \item $\Psi(00)=(0),\ \Psi(01)=(1)$.
      \item If $a=0\dots 0 12\in W^{n+1}$, we set $\Psi(a) = (0\dots011)\in \ma{N}^n$.
      \item If $a=A^{n-1} \cdots A^{2}(00)$, we set $L=B^{n-1} \cdots B^{2}(0)$.
      \item If $a=A^{n-1} \cdots A^{2}(01)$, we set $L=B^{n-1} \cdots B^{2}(1)$.
      \item If $a=A^{n-1} \cdots A^{k}(0\dots 012)$, we set $L=B^{n-1} \cdots B^{k}(0\dots 0 11)$.
   \end{enumerate}
   In particular, $|\mathcal N^n| = g_W(n+1)$ and
   \begin{equation*}
      \dim U(\mG_n)=\frac{(2^{n}+1)(2^{n-1}+1)}{3}\,.
   \end{equation*}
\end{theorem}
\begin{proof}
   From Theorem~\ref{thm:wordsconstuction} and the results in this section it is clear that $\Psi$ is a bijection, in particular it follows that
   $|\mathcal N^n| = g_W(n+1) = \frac{(2^{n}+1)(2^{n-1}+1)}{3}$.
   Therefore, the formula for $\dim U(\mG_n)$ follows from \eqref{eqr1} and \eqref{eq:clever}.
\end{proof}
\bigskip

\begin{example}
\label{ex:Psi}
   \begin{itemize}

      \item $\Psi(01010010) = 
      \begin{pmatrix}
	 1000000 \\ 0010000 \\ 0000010
      \end{pmatrix}
      $.
      \item 
      $\Psi(0123210) = 
      \begin{pmatrix}
	 101000\\ 010100 \\ 000010
      \end{pmatrix}$
      as the following diagram shows:

      \begin{equation*}
	 \begin{tikzcd}[column sep=1.5cm, row sep =1.2cm]
	    012         \arrow{r}[below]{\text{Case 5}}[above]{A^3_{3,2}} \arrow{d}[left]{\Psi}
	    & 0122      \arrow{r}[below]{\text{Case 6}}[above]{A^4_{4,3}}
	    & 01232     \arrow{r}[below]{\text{Case 2}}[above]{A^5_{6,1}}
	    & 012321    \arrow{r}[below]{\text{Case 1}}[above]{A^6_{7,0}}
	    & 0123210 \arrow{d}[right]{\Psi}
	    \\
	    (11) 
	    \arrow{r}[below]{\text{Case $5'$}}
	    & \begin{pmatrix} 101\\ 011
	    \end{pmatrix}
	    \arrow{r}[below]{\text{Case $6'a$}}
	    & \begin{pmatrix} 1010\\ 0101
	    \end{pmatrix}
	    \arrow{r}[below]{\text{Case $2'$}}
	    & \begin{pmatrix} 10100\\ 01010 \\ 00001
	    \end{pmatrix}
	    \arrow{r}[below]{\text{Case $1'$}}
	    & \begin{pmatrix} 101000\\ 010100 \\ 000010
	    \end{pmatrix}
	 \end{tikzcd}
      \end{equation*}


      \item 
      $\Psi(00122333) = 
      \begin{pmatrix} 0100000 \\ 0010100 \\ 0001100 \\ 0000011
      \end{pmatrix}$
      as the following diagram shows:

      \begin{equation*}
	 \begin{tikzcd}[column sep=1.5cm, row sep =1.2cm]
	    0012       \arrow{r}[below]{\text{Case 7}}[above]{A^4_{5,3}} \arrow{d}[left]{\Psi}
	    & 00123    \arrow{r}[below]{\text{Case 5}}[above]{A^5_{5,2}}
	    & 001223   \arrow{r}[below]{\text{Case 6}}[above]{A^6_{6,2}}
	    & 0012233  \arrow{r}[below]{\text{Case 6}}[above]{A^7_{7,3}}
	    & 00122333 \arrow{d}[right]{\Psi}
	    \\
	    (011) 
	    \arrow{r}[below]{\text{Case 7'}}
	    & \begin{pmatrix} 0100 \\ 0011
	    \end{pmatrix}
	    \arrow{r}[below]{\text{Case 5'}}
	    & \begin{pmatrix} 01000 \\ 00101 \\ 00011
	    \end{pmatrix}
	    \arrow{r}[below]{\text{Case 6'a}}
	    & \begin{pmatrix} 010000 \\ 001010 \\ 000101
	    \end{pmatrix}
	    \arrow{r}[below]{\text{Case 6'b}}
	    & \begin{pmatrix} 0100000 \\ 0010100 \\ 0001100 \\ 0000011
	    \end{pmatrix}
	 \end{tikzcd}
      \end{equation*}

   \end{itemize}
   \bigskip

   We give some examples of preimages of vector spaces (they are the spaces listed in \cite{li} on page~105).
   \begin{itemize}
      \item  
      $\Psi^{-1}
      \left(
      \begin{pmatrix}
	 10001 \\ 01000 \\ 00011
      \end{pmatrix}
      \right)
      = 011022$
      because

      \begin{equation*}
	 \begin{tikzcd}[column sep=1.5cm, row sep =1.2cm]
	    012      \arrow{r}[below]{\text{Case 4}}[above]{A^3_{3,1}} \arrow{d}[left]{\Psi}
	    & 0112   \arrow{r}[below]{\text{Case 3}}[above]{A^4_{4,0}}
	    & 01102  \arrow{r}[below]{\text{Case 5}}[above]{A^5_{5,2}}
	    & 011022 \arrow{d}[right]{\Psi}
	    \\
	    (11) 
	    \arrow{r}[below]{\text{Case $4'$}}
	    & \begin{pmatrix} 101 \\ 010
	    \end{pmatrix}
	    \arrow{r}[below]{\text{Case $3'$}}
	    & \begin{pmatrix} 1001 \\ 0100
	    \end{pmatrix}
	    \arrow{r}[below]{\text{Case $5'$}}
	    & \begin{pmatrix} 10001 \\ 01000 \\ 00011
	    \end{pmatrix}
	 \end{tikzcd}
      \end{equation*}

      \item 
      $\Psi^{-1}
      \left(
      \begin{pmatrix}
	 10010 \\ 01001 \\ 00101
      \end{pmatrix}
      \right)
      = 012232$
      because

      \begin{equation*}
	 \begin{tikzcd}[column sep=1.5cm, row sep =1.2cm]
	    012      \arrow{r}[below]{\text{Case 5}}[above]{A^3_{3,2}} \arrow{d}[left]{\Psi}
	    & 0122   \arrow{r}[below]{\text{Case 5}}[above]{A^4_{4,2}}
	    & 01222  \arrow{r}[below]{\text{Case 6}}[above]{A^5_{5,3}}
	    & 012232 \arrow{d}[right]{\Psi}
	    \\
	    (11) 
	    \arrow{r}[below]{\text{Case $5'$}}
	    & \begin{pmatrix} 101 \\ 011
	    \end{pmatrix}
	    \arrow{r}[below]{\text{Case $5'$}}
	    & \begin{pmatrix} 1001 \\ 0101 \\ 0011
	    \end{pmatrix}
	    \arrow{r}[below]{\text{Case $6'a$}}
	    & \begin{pmatrix} 10010 \\ 01001 \\ 00101
	    \end{pmatrix}
	 \end{tikzcd}
      \end{equation*}

      \item  
      $\Psi^{-1}
      \left(
      \begin{pmatrix}
	 1000010 \\ 0100001 \\ 0010001 \\ 0001000 \\ 0000101
      \end{pmatrix}
      \right)
      = 01221232$
      because

      \begin{equation*}
	 \hspace*{-1cm}
	 \begin{tikzcd}[column sep=1.2cm, row sep =1.2cm]
	    012        \arrow{r}[below]{\text{Case 5}}[above]{A^3_{3,2}} \arrow{d}[left]{\Psi}
	    & 0122     \arrow{r}[below]{\text{Case 5}}[above]{A^4_{4,2}}
	    & 01222    \arrow{r}[below]{\text{Case 4}}[above]{A^5_{5,1}}
	    & 012212   \arrow{r}[below]{\text{Case 5}}[above]{A^6_{6,2}}
	    & 0122122  \arrow{r}[below]{\text{Case 6}}[above]{A^7_{7,3}}
	    & 01221232 \arrow{d}[right]{\Psi}
	    \\
	    (11) 
	    \arrow{r}[below]{\text{Case $5'$}}
	    & \begin{pmatrix} 101 \\ 011
	    \end{pmatrix}
	    \arrow{r}[below]{\text{\!\!\!Case $5'$}}
	    & \begin{pmatrix} 1001 \\ 0101 \\ 0011
	    \end{pmatrix}
	    \arrow{r}[below]{\text{Case $4'$}}
	    & \begin{pmatrix} 10001 \\ 01001 \\ 00101 \\ 00010
	    \end{pmatrix}
	    \arrow{r}[below]{\text{Case $5'$}}
	    & \begin{pmatrix} 100001 \\ 010001 \\ 001001 \\ 000100 \\ 000011
	    \end{pmatrix}
	    \arrow{r}[below]{\text{Case $6'a$}}
	    & \begin{pmatrix} 1000010 \\ 0100001 \\ 0010001 \\ 0001000 \\ 0000101 
	    \end{pmatrix}
	 \end{tikzcd}
      \end{equation*}

   \end{itemize}

\end{example}


\appendix

\section{Symplectic Dual Polar Space for $n=1,2,3$}
\label{appendix:CR}
Recall that $\mathcal G_n$ is the symplectic dual polar space defined in Section~\ref{sec:symplectic}.
It consists of $\prod_{k=1}^{n}(2^k+1)$ points (Corollary~\ref{cor:NN}) and $\frac{1}{3} (2^n-1)\prod_{k=1}^{n}(2^k+1)$ lines (Corollary~\ref{cor:NHP}).
Each line contains exactly three points by Corollary~\ref{cor:tres} 
and through each point pass exactly $2^n-1$ lines because any $n$-dimensional $\FF_2$-vector space contains exactly $2^n-1$ different $(n-1)$-dimensional subspaces.

\subsection{$n=1$} 
$\mathcal G_1$ consists of one line with exactly three points.

\subsection{$n=2$} 
For $\mathcal G_2$ we label the 15 points by 
\begin{alignat*}{6} 
   0&\leftrightarrow\begin{pmatrix} 
      1 0 0 0\\ 0 1 0 0
   \end{pmatrix},\quad
   &1&\leftrightarrow\begin{pmatrix} 
      1 0 0 0\\ 0 1 0 1
   \end{pmatrix},\quad
   &2&\leftrightarrow\begin{pmatrix} 
      1 0 0 1\\ 0 1 1 0
   \end{pmatrix},\quad
   &3&\leftrightarrow\begin{pmatrix} 
      1 0 0 1 \\ 0 1 1 1
   \end{pmatrix},\quad
   &4&\leftrightarrow\begin{pmatrix} 
      1 0 1 0 \\ 0 1 0 0
   \end{pmatrix},
   \\[1ex]
   5&\leftrightarrow\begin{pmatrix} 
      1 0 1 0 \\ 0 1 0 1
   \end{pmatrix},
   &6&\leftrightarrow\begin{pmatrix} 
      1 0 1 1\\ 0 1 1 0
   \end{pmatrix},
   &7&\leftrightarrow\begin{pmatrix} 
      1 0 1 1\\ 0 1 1 1
   \end{pmatrix},
   &8&\leftrightarrow\begin{pmatrix} 
      1 1 0 0\\ 0 0 1 1
   \end{pmatrix},
   &9&\leftrightarrow\begin{pmatrix} 
      1 1 0 1\\ 0 0 1 1
   \end{pmatrix},
   \\[1ex]
   10&\leftrightarrow\begin{pmatrix} 
      1 0 0 0\\ 0 0 0 1
   \end{pmatrix},
   &11&\leftrightarrow\begin{pmatrix}
      1 0 1 0\\ 0 0 0 1
   \end{pmatrix},
   &12&\leftrightarrow\begin{pmatrix} 
      0 1 0 0\\ 0 0 1 0
   \end{pmatrix},
   &13&\leftrightarrow\begin{pmatrix} 
      0 1 0 1\\ 0 0 1 0
   \end{pmatrix},
   &14&\leftrightarrow\begin{pmatrix} 
      0 0 1 0\\ 0 0 0 1
   \end{pmatrix}.
\end{alignat*}
     We used SageMath \cite{sag} to write down all lines as triples of their points:
\begin{center}   
   \begin{tabular}[t]{lllll}
      (0, 1, 10), & (6, 7, 11), &  (3, 4, 9), & (1, 5, 13), & (12, 13, 14), \\
      (2, 3, 10), &  (0, 7, 8), &  (2, 5, 8), & (2, 6, 12), & (8, 9, 14), \\
      (4, 5, 11), &  (1, 6, 9), & (0, 4, 12), & (3, 7, 13), & (10, 11, 14).
   \end{tabular}
\end{center}
We fix the vertex $0$ and construct the corresponding subgraphs $\Gamma_0,\Gamma_1$ and $\Gamma_2$. 

\begin{itemize}
\item $\Gamma_0$ consists only of the vertex $0$.
\item $\Gamma_1$ consists of three connected components, each of which contains two vertices and one edge.
\item $\Gamma_2$ consists of one connected component and eight points which form a cube.
\end{itemize}
The subgraphs $\Gamma_1$ and $\Gamma_2$ are shown in Figure \ref{fif2}.
\begin{figure}[h!]
   \begin{center}
      \begin{tikzpicture}[scale=.5] 
	 \begin{scope}
	    \tikzset{transform shape};
	    \tikzset{esquina/.style = {draw, circle, minimum size=#1, inner sep=0pt, outer sep=0pt},
	    esquina/.default = 18pt  
	    }
	    \tikzset{mycube/.pic={%
	    \node[esquina, name=a1] at (0,0) {14}; 
	    \node[esquina, name=a2] at (4,0) {9};  
	    \node[esquina, name=a3] at (4,4) {3};  
	    \node[esquina, name=a4] at (0,4) {13}; 
	    \node[esquina, name=a5] at (1,2) {11}; 
	    \node[esquina, name=a6] at (5,2) {6};  
	    \node[esquina, name=a7] at (5,6) {2};  
	    \node[esquina, name=a8] at (1,6) {5};  

	    \draw[very thick] (a1)--(a2);
	    \draw[very thick] (a2)--(a3);
	    \draw[very thick] (a3)--(a4);
	    \draw[very thick] (a4)--(a1);

	    \draw[very thick] (a4)--(a8);
	    \draw[very thick] (a3)--(a7);
	    \draw[very thick] (a2)--(a6);
	    \draw[very thick] (a6)--(a7);
	    \draw[very thick] (a7)--(a8);

	    \draw[very thick, dashed] (a5)--(a1);
	    \draw[very thick, dashed] (a5)--(a6);
	    \draw[very thick, dashed] (a5)--(a8);
	    }}

	    \node [esquina, name=g1] at (0,0) {1};
	    \node [esquina, name=g2] at (0,5) {10};
	    \node [esquina, name=g3] at (3,0) {4};
	    \node [esquina, name=g4] at (3,5) {12};
	    \node [esquina, name=g5] at (6,0) {7};
	    \node [esquina, name=g6] at (6,5) {8};
	    \draw[very thick] (g1)--(g2);
	    \draw[very thick] (g3)--(g4);
	    \draw[very thick] (g5)--(g6);

	    \pic at (15,0) {mycube};
	 \end{scope}
	 \node [anchor=west] at (-2, -3) {%
	 \parbox{.3\textwidth}{
	 Subgraph $\Gamma_1$ of $\mathcal G_2$ consisting of three connected components.}};
	 \node [anchor=west] at (12, -3) {%
	 \parbox{.3\textwidth}{
	 Subgraph $\Gamma_2$ of $\mathcal G_2$ consisting of one connected component.}};

      \end{tikzpicture} 
   \end{center}
   \caption{Induced subgraphs $\Gamma_1$ and $\Gamma_2$ of $\mathcal G_2$.}
   \label{fif2}
\end{figure}

\subsection{$n=2$} 
The dual polar space $\mathcal G_3$ consists of 135 points and 315 lines.
Each line contains exactly three points and through each point pass exactly seven lines.
The subgraphs $\Gamma_0,\Gamma_1$, $\Gamma_2$ and $\Gamma_3$ have the following description:
\begin{itemize}
    \item $\Gamma_0$ consists only of the vertex $0$.
\item $\Gamma_1$ has seven connected components each of which consists of two vertices and one edge.
\item $\Gamma_2$ consists of seven components each one with the form of a cube.
\item $\Gamma_3$ is connected and consists of $64$ vertices and $224$ edges.
\end{itemize}

\section{Construction of $W^{n+1}$ from $W^n$}
\label{appendix:construction}
Let $n\ge 2$.
In this appendix we show how $W^{n+1}$ is constructed from $W^n$.
Recall that if $a=a_1a_2\dots a_{n-1} a_n\in W^n$, then we can do the following (cf. Remark~\ref{rem:cases}):
\begin{itemize}
   \item If $a_n\in\{0,1\}$: append $0$ or $1$. We obtain a word in Case 1 or 2.
   \item If $a_n\in\{2,3\}$: insert $0$, $1$ or $2$ before $a_n$. We obtain a word in Case 3, 4 or 5.
   \item If $a_n\in\{2,3\}$ and it is possible to insert $3$ before $a_n$, we do so. We obtain a word in Case 6.
   \item If $a_n\in\{2,3\}$ and it is not possible to insert $3$ before $a_n$, then necessarily $a_n=2$ and $a_j\in\{0,1\}$ for $1\le j \le n-1$. We append $3$ and obtain a word in Case 7.
\end{itemize}
Finally, we have to add the word $a=0\dots 012$ from Case~8.
In this way we obtain all possible words of $W^{n+1}$.
\medskip

\newcommand{\mc}[2]{\multicolumn{#1}{c}{#2}}
\newcommand{\C}[1]{\multicolumn{#1}{c}{}}
\newcommand{\Cl}[1]{\multicolumn{#1}{c|}{}}
\newcommand{\CA}[1]{\multicolumn{#1}{c|}{$\downarrow$}}

\begin{minipage}{10cm}
\rotatebox{-90}{
\small
\begin{tabular}[t]{|l|l|*{15}{l|}l|}
   \cline{1-17}
   \bf n=4 & last & \multicolumn{10}{c|}{$0$ or $1$} &\multicolumn{5}{c|}{$2$ or $3$} & \C{1} \\
   \cline{1-17}
   && 0000 & 0010 & 0100 & 0110 & 0120 &
      0001 & 0011 & 0101 & 0111 & 0121 &
      0102 & 0112 & 0122 & 0123 & 0012 &
   \C{1} \\
   \cline{1-17}
   \multicolumn{2}{|c|}{} & 
   \CA{1} & \CA{1}& \CA{1}& \CA{1} & \CA{1} & 
   \CA{1} & \CA{1}& \CA{1}& \CA{1} & \CA{1} & 
   \CA{1} & \CA{1}& \CA{1}& \CA{1} & \CA{1} & 
   \C{1} \\
   \cline{1-4}
   \hline
   \bf n=5 
   & C1 & 
   0000\m{0} & 0010\m{0} & 0100\m{0} & 0110\m{0} & 0120\m{0} &
   0001\m{0} & 0011\m{0} & 0101\m{0} & 0111\m{0} & 0121\m{0} &
   0102\m{0} & 0112\m{0} & 0122\m{0} & 0123\m{0} & 0012\m{0} &
   \Cl{1} \\ \hline
   & C2 & 
   0000\m{1} & 0010\m{1} & 0100\m{1} & 0110\m{1} & 0120\m{1} &
   0001\m{1} & 0011\m{1} & 0101\m{1} & 0111\m{1} & 0121\m{1} &
   0102\m{1} & 0112\m{1} & 0122\m{1} & 0123\m{1} & 0012\m{1} &
   \Cl{1} \\ \cline{2-18}
   & C3 & & & & & & & & & & & 010\m{0}2 & 011\m{0}2 & 012\m{0}2 & 012\m{0}3 & 001\m{0}2 &\Cl{1}\\ \cline{2-18}
   & C4 & & & & & & & & & & & 010\m{1}2 & 011\m{1}2 & 012\m{1}2 & 012\m{1}3 & 001\m{1}2 &\Cl{1}\\ \cline{2-18}
   & C5 & & & & & & & & & & & 010\m{2}2 & 011\m{2}2 & 012\m{2}2 & 012\m{2}3 & 001\m{2}2 &\Cl{1}\\ \cline{2-18}
   & C6 & & & & & & & & & & &           &           & 012\m{3}2 & 012\m{3}3 &           &\Cl{1}\\ \cline{2-18}
   & C7 & & & & & & & & & & & 0102\m{3} & 0112\m{3} &           &           & 0012\m{3} &\Cl{1}\\ \cline{2-18}
   & C8 & & & & & & & & & & & & & & & & 00012 \\
   \hline
\end{tabular}}
\end{minipage}
\               \
\begin{minipage}{5cm}
\rotatebox{-90}{
\begin{tabular}[t]{|l|l|l|l|l|}
   \cline{1-4}
   \bf n=2 & last letter & \multicolumn{2}{c|}{$0$ or $1$} & \C{1} \\
   \cline{1-4}
   && 00 & 01 & \C{1} \\
   \cline{1-4}
   \multicolumn{2}{|c|}{} & \CA{1} & \CA{1} & \C{1} \\
   \cline{1-4}
   \hline
   \bf n=3 & Case 1 & 00\m{0} & 01\m{0} & \Cl{1} \\
   \hline
   & Case 2 & 00\m{1} & 01\m{1} & \Cl{1} \\
   \cline{2-5}
   & Case 8 &  &  & 012 \\
   \hline
\end{tabular}}
\rotatebox{-90}{
\begin{tabular}[t]{|l|l|*{5}{l|}l|}
   \cline{1-7}
   \bf n=3 & last letter & \multicolumn{4}{c|}{$0$ or $1$} & $2$\! or\! $3$ & \C{1} \\
   \cline{1-7}
   && 000 & 001 & 010 & 011 & 012 & \C{1} \\
   \cline{1-7}
   \multicolumn{2}{|c|}{} & \CA{1} & \CA{1}& \CA{1}& \CA{1} & \CA{1} & \C{1} \\
   \cline{1-4}
   \hline
   \bf n=4 
   & Case 1 & 000\m{0} & 001\m{0} & 010\m{0} & 011\m{0} & 012\m{0} & \Cl{1} \\ \hline
   & Case 2 & 000\m{1} & 001\m{1} & 010\m{1} & 011\m{1} & 012\m{1} & \Cl{1} \\ \cline{2-8}
   & Case 3 &          &          &          &          & 01\m{0}2 & \Cl{1} \\ \cline{2-8}
   & Case 4 &          &          &          &          & 01\m{1}2 & \Cl{1} \\ \cline{2-8}
   & Case 5 &          &          &          &          & 01\m{2}2 & \Cl{1} \\ \cline{2-8}
   & Case 6 &          &          &          &          &          & \Cl{1} \\ \cline{2-8}
   & Case 7 &          &          &          &          & 012\m{3} & \Cl{1} \\ \cline{2-8}
   & Case 8 & & & & &  & 0012 \\
   \hline
\end{tabular}}
\bigskip
\end{minipage}


\section{Classification of words in $W^n$ according to the Cases 1 -- 8}
\label{apenC}
Let $n\ge 2$ and $a=a_1a_2\dots a_{n-1} a_n\in W^n$.
As before, we set $E^n_{n-1}(a) = a_1a_2\dots a_{n-2}a_n$ which is obtained from $a$ by erasing its second to last letter.
Recall that the word $a$ belongs to
\begin{itemize}
   \addtolength{\itemsep}{-.5\baselineskip}
   \item {\bf Case 1} if $a_n=0$;
   \item {\bf Case 2} if $a_n=1$;
   \item {\bf Case 3} if $a_n\in\{2,3\}$ and $a_{n-1}=0$ (then automatically $E^n_{n-1}(a)\in W^{n-1}$);
   \item {\bf Case 4} if $a_n\in\{2,3\}$, $a_{n-1}=1$ and $E^n_{n-1}(a)\in W^{n-1}$;
   \item {\bf Case 5} if $a_n\in\{2,3\}$, $a_{n-1}=2$ and $E^n_{n-1}(a)\in W^{n-1}$;
   \item {\bf Case 6} if $a_n\in\{2,3\}$ and $a_{n-1}=3$ (then automatically $E^n_{n-1}(a)\in W^{n-1}$);
   \item {\bf Case 7} if $a_n\in\{2,3\}$, $a_{n-1}=2$ and $E^n_{n-1}(a)\not\in W^{n-1}$\\
   \phantom{\bf Case 7 }(equivalently: if $a_{n-1}a_n=23$ and $a_j\in\{0,1\}$ for $1\le j\le n-2$);
   \item {\bf Case 8} if $a_n\in\{2,3\}$, $a_{n-1}=1$ and $E^n_{n-1}(a)\not\in W^{n-1}$\\
   \phantom{\bf Case 8 }(equivalently: if $a = 0\dots 0 12$).
\end{itemize}

So we obtain for $n=1,2,3,4,5$:
\bigskip

\begin{minipage}[t]{.1\textwidth}
   \fbox{n=1}
\end{minipage}
\begin{minipage}[t]{.7\textwidth}
   \vspace*{-1\baselineskip}
   \begin{tabular}[t]{| r | r | r | r | r | r | r | r |}
      \hline
      Case 1 & Case 2 & Case 3 & Case 4 & Case 5 & Case 6 & Case 7 & Case 8\\ \hline
      0 &  & & & & & & \\
      \hline
   \end{tabular}
\end{minipage}
\bigskip

\begin{minipage}[t]{.1\textwidth}
   \fbox{n=2}
\end{minipage}
\begin{minipage}[t]{.7\textwidth}
   \vspace*{-1\baselineskip}
   \begin{tabular}[t]{| r | r | r | r | r | r | r | r |}
      \hline
      Case 1 & Case 2 & Case 3 & Case 4 & Case 5 & Case 6 & Case 7 & Case 8\\ \hline
      0\m{0} & 0\m{1} & & & & & & \\
      \hline
   \end{tabular}
\end{minipage}
\bigskip

\begin{minipage}[t]{.1\textwidth}
   \fbox{n=3}
\end{minipage}
\begin{minipage}[t]{.7\textwidth}
   \vspace*{-1\baselineskip}
   \begin{tabular}[t]{| r | r | r | r | r | r | r | r |}
      \hline
      Case 1  & Case 2  & Case 3 & Case 4 & Case 5 & Case 6 & Case 7 & Case 8\\ \hline
      00\m{0} & 00\m{1} & & & & & & 012 \\
      01\m{0} & 01\m{1} & & & & & &     \\
      \hline
   \end{tabular}
\end{minipage}
\bigskip

\begin{minipage}[t]{.1\textwidth}
   \fbox{n=4}
\end{minipage}
\begin{minipage}[t]{.7\textwidth}
   \vspace*{-1\baselineskip}
   \begin{tabular}[t]{| r | r | r | r | r | r | r | r |}
      \hline
      Case 1   & Case 2   & Case 3   & Case 4   & Case 5   & Case 6 & Case 7   & Case 8\\ \hline
      000\m{0} & 000\m{1} & 01\m{0}2 & 01\m{1}2 & 01\m{2}2 &        & 012\m{3} & 0012   \\
      001\m{0} & 001\m{1} &          &          &          &        &          &        \\
      010\m{0} & 010\m{1} &          &          &          &        &          &        \\
      011\m{0} & 011\m{1} &          &          &          &        &          &        \\
      012\m{0} & 012\m{1} &          &          &          &        &          &        \\
      \hline
   \end{tabular}
\end{minipage}
\bigskip

\begin{minipage}[t]{.1\textwidth}
   \fbox{n=5}
\end{minipage}
\begin{minipage}[t]{.7\textwidth}
   \vspace*{-1\baselineskip}
   \begin{tabular}[t]{| r | r | r | r | r | r | r | r |}
      \hline
      Case 1    & Case 2    & Case 3    & Case 4    & Case 5    & Case 6    & Case 7    & Case 8\\ \hline
      0000\m{0} & 0000\m{1} & 010\m{0}2 & 010\m{1}2 & 010\m{2}2 &           & 0102\m{3} & 00012   \\
      0010\m{0} & 0010\m{1} & 011\m{0}2 & 011\m{1}2 & 011\m{2}2 &           & 0112\m{3} &         \\
      0100\m{0} & 0100\m{1} & 012\m{0}2 & 012\m{1}2 & 012\m{2}2 & 012\m{3}2 &           &         \\
      0110\m{0} & 0110\m{1} & 012\m{0}3 & 012\m{1}3 & 012\m{2}3 & 012\m{3}3 &           &         \\
      0120\m{0} & 0120\m{1} & 001\m{0}2 & 001\m{1}2 & 001\m{2}2 &           & 0012\m{3} &         \\
      0001\m{0} & 0001\m{1} &           &           &           &           &           &         \\
      0011\m{0} & 0011\m{1} &           &           &           &           &           &         \\
      0101\m{0} & 0101\m{1} &           &           &           &           &           &         \\
      0111\m{0} & 0111\m{1} &           &           &           &           &           &         \\
      0121\m{0} & 0121\m{1} &           &           &           &           &           &         \\
      0102\m{0} & 0102\m{1} &           &           &           &           &           &         \\
      0112\m{0} & 0112\m{1} &           &           &           &           &           &         \\
      0122\m{0} & 0122\m{1} &           &           &           &           &           &         \\
      0123\m{0} & 0123\m{1} &           &           &           &           &           &         \\
      0012\m{0} & 0012\m{1} &           &           &           &           &           &         \\
      \hline
   \end{tabular}
\end{minipage}
\newpage
{\bf Acknowledgemnt.}
C. Segovia is supported by catedras CONACYT and Proyecto CONACYT ciencias b\'asicas 2016, No. 284621.

\bibliographystyle{alpha}
\bibliography{lit-SeWi-2017}

\end{document}

%% file: CR.tex
\begin{tikzpicture}[scale=.6] 
      \node[draw,circle, violet, name=a] at (0*72:1) {};  
      \node[draw,circle, violet, name=b] at (1*72:1) {};  
      \node[draw,circle, violet, name=c] at (2*72:1) {};  
      \node[draw,circle, violet, name=d] at (3*72:1) {};  
      \node[draw,circle, violet, name=e] at (4*72:1) {};  

      \node[draw,circle, violet, name=aaa] at (0*72:5.7) {};   
      \node[draw,circle, violet, name=bbb] at (1*72:5.7) {};   
      \node[draw,circle, violet, name=ccc] at (2*72:5.7) {};   
      \node[draw,circle, violet, name=ddd] at (3*72:5.7) {};   
      \node[draw,circle, violet, name=eee] at (4*72:5.7) {};   


      \path[green, name path=aaab] (aaa)--($(b)!-8cm!(aaa)$);
      \path[red,   name path=ac] (a)--($(c)!-6cm!(a)$);
      \path [name intersections={of={aaab and ac}}];
      \coordinate (rrr) at (intersection-1);
      \node[draw, circle, violet, name=R] at (rrr) {};

      \path[green, name path=bbbc] (bbb)--($(c)!-8cm!(bbb)$);
      \path[red,   name path=bd] (b)--($(d)!-6cm!(b)$);
      \path [name intersections={of={bbbc and bd}}];
      \coordinate (sss) at (intersection-1);
      \node[draw, circle, violet, name=S] at (sss) {};

      \path[green, name path=cccd] (ccc)--($(d)!-8cm!(ccc)$);
      \path[red,   name path=ce] (c)--($(e)!-6cm!(c)$);
      \path [name intersections={of={cccd and ce}}];
      \coordinate (ttt) at (intersection-1);
      \node[draw, circle, violet, name=T] at (ttt) {};

      \path[green, name path=ddde] (ddd)--($(e)!-8cm!(ddd)$);
      \path[red,   name path=da] (d)--($(a)!-6cm!(d)$);
      \path [name intersections={of={ddde and da}}];
      \coordinate (uuu) at (intersection-1);
      \node[draw, circle, violet, name=U] at (uuu) {};

      \path[green, name path=eeea] (eee)--($(a)!-8cm!(eee)$);
      \path[red,   name path=eb] (e)--($(b)!-6cm!(e)$);
      \path [name intersections={of={eeea and eb}}];
      \coordinate (vvv) at (intersection-1);
      \node[draw, circle, violet, name=V] at (vvv) {};

      \draw[thick] (a)--(c); \draw[thick](c)--(R);
      \draw[thick] (b)--(d); \draw[thick](d)--(S);
      \draw[thick] (c)--(e); \draw[thick](e)--(T);
      \draw[thick] (d)--(a); \draw[thick](a)--(U);
      \draw[thick] (e)--(b); \draw[thick](b)--(V);

      \draw[thick] (aaa)--(T); \draw[thick](T)--(eee);
      \draw[thick] (bbb)--(U); \draw[thick](U)--(aaa);
      \draw[thick] (ccc)--(V); \draw[thick](V)--(bbb);
      \draw[thick] (ddd)--(R); \draw[thick](R)--(ccc);
      \draw[thick] (eee)--(S); \draw[thick](S)--(ddd);

      \draw[thick] (aaa)--(b); \draw[thick](b)--(R);
      \draw[thick] (bbb)--(c); \draw[thick](c)--(S);
      \draw[thick] (ccc)--(d); \draw[thick](d)--(T);
      \draw[thick] (ddd)--(e); \draw[thick](e)--(U);
      \draw[thick] (eee)--(a); \draw[thick](a)--(V);

   \end{tikzpicture}